\documentclass[12pt]{article} 
\usepackage{amssymb}
\usepackage{graphicx}
\def\bfc{{\bf C}}
\def\bfr{{\bf R}}

\def\Bar{\overline}

\def\dim{{\rm dim}}
\def\diam{{\rm diam}}

\def\grad{{\rm grad}}
\def\hess{{\rm Hess}}

\def\ident{\equiv}

\def\intersect{\bigcap}
\def\iso{\cong}

\def\na{\nabla}

\def\span{{\rm span}} 

\def\tensor{\otimes}

\def\Union{\bigcup}

\def\be{\begin{equation}}
\def\ee{\end{equation}}
\def\bearray{\begin{eqnarray}}
\def\eearray{\end{eqnarray}}
\def\bestar{\begin{eqnarray*}}
\def\eestar{\end{eqnarray*}}
\def\ben{\begin{displaymath}}
\def\een{\end{displaymath}}

\newtheorem{theorem}{Theorem}[section]
\newtheorem{thm}[theorem]{Theorem}

\newtheorem{prop}[theorem]{Proposition}
\newtheorem{lemma}[theorem]{Lemma}
\newtheorem{corollary}[theorem]{Corollary}
\newtheorem{cor}[theorem]{Corollary}
\newtheorem{defn}[theorem]{Definition}

\newtheorem{remark}[theorem]{Remark}

\newenvironment{example}%
{\addtocounter{theorem}{1}%
\noindent {\bf
Example \arabic{section}.\arabic{theorem} }}{\\[2 ex]}

\def\qed{{\hfill\vrule height10pt width10pt}\\ \vspace{3 ex}}
\def\qedns{{\hfill\vrule height10pt width10pt} \\ \vspace{0.5 ex}
}

\def\pf{\ni{\bf Proof}: }

\def\d{\delta}
\def\e{\epsilon}
\def\g{\gamma}
\def\i{\iota}
\def\k{\kappa}
\def\l{\lambda}

\def\ni{\noindent}
\def\ss{\vspace{.1in}}
\def\bs{\vspace{.2in}}
\def\ms{\vspace{.15in}}
\def\bsn{\vspace{.2in}\noindent}

\def\ssn{\vspace{.1in}\noindent}

\def\hor{{\rm hor}}
\def\vert{{\rm vert}}
\def\tg{\tilde{\g}}

\def\aku{|K|}
\def\P{{\cal P}}
\def\It{{\rm It}}

\def\rinj{r_{\rm inj}}
\def\rreg{r_{\rm reg}}
\def\rcvx{r_{\rm cvx}}
\def\rregcvx{r_{\rm regcvx}}
\def\circumrad{{\rm circumrad}}

\def\hess{{\rm Hess}}
\def\hessp{{\rm Hess}'}
\def\ch{{\bf c}}
\def\sh{{\bf s}}
\def\psimax{\psi_{\rm max}}

\def\tq{{\tilde{Q}}}
\def\tQ{{\tilde{Q}}}

\def\hull{{\rm hull}}
\def\ohull{{\rm ohull}}
\def\Star{{\rm star}}
\def\regstar{{\rm regstar}}

\begin{document}


\setcounter{page}{0}

\date{Revised version: July 2, 2003}
\title{\bf 
Newton's method, zeroes of vector fields, and the 
Riemannian center of mass}
\author{David Groisser
\\ {\small Department of Mathematics}
\\ {\small University of Florida}
\\ {\small Gainesville FL 32611--8105}
\\ {\small USA}
\\ {\small groisser@math.ufl.edu}}
\maketitle

\begin{abstract}

We present an iterative technique for finding zeroes of vector fields
on Riemannian manifolds. As a special case we obtain a ``nonlinear
averaging algorithm'' that computes the centroid of a mass
distribution $\mu$ supported in a set of small enough diameter $D$ in
a Riemannian manifold $M$.  We estimate the convergence rate of our
general algorithm and the more special Riemannian averaging algorithm.
The algorithm is also used to provide a constructive proof of
Karcher's theorem on the existence and local uniqueness of the center
of mass, under a somewhat stronger requirement than Karcher's on $D$.
Another corollary of our results is a proof of
convergence, for a fairly large open set of initial conditions, of the
``GPA algorithm'' used in statistics to average points in a
shape-space, and a quantitative explanation of why the GPA
algorithm converges rapidly in practice; see \cite{groi2}.

We also show that a mass distribution in $M$ with support $Q$ has a
unique center of mass in a (suitably defined) convex hull of $Q$.

\end{abstract}


\vspace{.2in}
   

\newpage
\section{Introduction.}

In this article we present an iterative technique for finding zeroes
of vector fields on Riemannian manifolds, and apply this technique to
the averaging of a mass distribution with support contained in a
sufficiently small ball in a Riemannian manifold.  Our approach
provides a new and constructive proof of Karcher's theorem on the
existence and uniqueness of the center of mass, under a somewhat
stronger requirement on the radius of the supporting ball than was
used in \cite{karcher}.

This study was originally motivated by curiosity about a method (the
``GPA algorithm'') used in statistics to find the average, suitably
defined, of a sample of shapes.  In many areas of image analysis,
particularly in biological applications such as cardiography
(cf. \cite{sbsb}) and maps of the brain (cf. \cite{bookstein4}) this
average is the starting point for understanding ``normal'' shapes and
deviations from the norm.  In practical applications the averaging
algorithm tends to converge remarkably quickly, often stabilizing to
desired precision after two or three iterations
(cf. \cite{bookstein4}, Figure 5 (p. 22), or \cite{goodall} Table 3
(p. 307)).  The initial purpose of our study was to understand the
geometry underlying this algorithm and, in quantitative terms, why the
convergence in practical applications is so rapid.  In exploring this
the author found that the GPA algorithm has a more general
interpretation on Riemannian manifolds, generalizing to a technique
for finding local zeroes of a vector field. The technique is an
iterative algorithm that we show is closely related to Newton's method
and mimics the contracting-mapping proof of the Inverse Function
Theorem.

As a special case of this technique, we obtain a general Riemannian
averaging algorithm. The vector field used in this algorithm has a
unique local zero, assuming the diameter $D$ of the support of
distribution being averaged is not too large, and is ``almost linear''
near this zero if $D$ is small, explaining the rapid convergence.
This zero is exactly the Riemannian center of mass of the distribution
being averaged.  In sections \ref{construct} and
\ref{conv_rate} of this paper we quantify ``not too large'' and
``small'', giving sufficient conditions for convergence of the
algorithm and estimating the convergence rate.

The Riemannian averaging algorithm can in principle be applied to any
``nonlinear averaging'' problem in which the objects being averaged
are parametrized by a Riemannian manifold, and is easily implemented
in spaces for which the exponential map and its inverse are explicitly
known (e.g. Riemannian submersions from spheres, and certain
homogeneous spaces with invariant metrics).  This is exactly the
situation for the shape-averaging problem. The (Euclidean) {\em shape
space} $\Sigma^k_n$ is the space of configurations of $k$
non-identical labeled points in $\bfr^n$, modulo equivalence under
translations, rotations, and dilations (rescalings) in $\bfr^n$;
sometimes one also allows reflections. The {\em size-and-shape space}
$\tilde{\Sigma}^k_n$ is defined similarly, but one does not mod out by
rescalings.  These spaces can naturally be given the structure of
manifolds with singularities, with natural Riemannian metrics on their
smooth parts (\cite{dkendall, carne, kendall_le}).  Averaging
(sizes-and-) shapes can be viewed as averaging certain mass
distributions on (size-and-) shape spaces, namely finite lists of
points with normalized counting measure.  In the probability and
statistics literature there is a commonly accepted definition of mean
size-and-shape, the Procrustean mean size-and-shape, but several
possible definitions of mean shape (see \cite{kent1} p. 292 and
\cite{kent2}), the most common of which may be the Procrustean mean shape
used in \cite{Le1}.  However, while the Procrustean mean {\em
size-and-shape} as defined in the probability and statistics
literature agrees with the Riemannian center of mass, the Procrustean
mean {\em shape} does not.

The GPA (Generalized Procrustes Analysis) algorithm as described in
\cite{Le1} lives intrinsically on size-and-shape space; call this
algorithm GPA-SS.  To obtain from this an algorithm that averages
shapes, one first embeds shape-space $\Sigma^k_n$ into size-and-shape
space $\tilde{\Sigma}^k_n$ in a standard way, carrying the list $Q$ of
shapes to be averaged to a list $\i(Q)$ of sizes-and-shapes. One then
produces a sequence of in $\tilde{\Sigma}^k_n$ by applying the GPA-SS
algorithm to $\i(Q)$. Finally one projects the limit (if there is one)
back onto shape space. Call this set of steps GPA-S. Le proves in
\cite{Le1} that if the shapes in $Q$ 
are not too far apart in $\Sigma^k_n$, and if the sequence in
$\tilde{\Sigma}^k_n$ converges, then the limit in $\tilde{\Sigma}^k_n$
is the Procrustean mean size-and-shape of the list $\i(Q)$.  It is not
hard to show that this projection of the Procrustean mean
size-and-shape is exactly the Procrustean mean shape (\cite{Le1},
p. 54), so that GPA-S computes the Procrustean mean shape.

Although the literature contains many discussions of the GPA-SS and
other GPA-derived algorithms, at the time this paper was first
completed \cite{groi1} the literature contained no theorems giving
sufficient conditions for any of these algorithms to converge.
However, as we show in \cite{groi2}, the GPA-SS algorithm is {\em
exactly} our Riemannian averaging algorithm as applied to
size-and-shape space. Hence convergence of the GPA-SS algorithm, for
an explicitly describable open set of initial conditions, is an
immediate corollary of the Riemannian-averaging theorems in sections
\ref{construct} and \ref{conv_rate} of this paper. After \cite{groi1}
was written, \cite{Le3}, which contains some overlapping results,
appeared.

In the iterative part of the GPA-S algorithm, one can obtain a
sequence of points in shape space by projecting each point in the
GPA-SS sequence, rather than just the limit, back onto shape
space. (This sequence in shape space can also be described slightly
more intrinsically; see \cite{groi2}, where we discuss the application
of the results of this paper to Procrustean averaging in more detail.)
In this way one obtains an iterative algorithm GPA-S$'$ on shape space
itself. GPA-S$'$ does not coincide with the Riemannian averaging
algorithm on shape space---it cannot, since it converges (for suitable
initial conditions) to the Procrustean mean shape and not to the
Riemannian average. However, GPA-S$'$ is an algorithm of the more
general type also considered here, and therefore its convergence,
again for an explicitly describable open set of initial conditions,
follows directly from our more general theorems in section
\ref{zvf}, as well as from the fact that GPA-S converges. 

In this paper we also address the question of why the convergence of
the GPA algorithms is so rapid in practice. As has been noted by many
authors, the data sets averaged in practical applications tend to be
very concentrated sets in shape (or size-and-shape) space; their
diameter $D$ is very small compared with any length-scale derivable
from the geometry of shape (or size-and-shape) space.  Our theorems in
section \ref{conv_rate} show why, for small $D$, convergence is rapid.

To describe our results more concretely, we need some notation and
terminology:

\begin{defn}\label{defcont} {\rm 
Let $(A,d_A), (B,d_B)$ be metric spaces and let $\k\in[0,1).$ We 
call a map $F:A\to B$ a {\em contraction with constant $\k$} if
$d_B(F(x),F(y))\leq \k\ d_A(x,y)$ for all $x,y\in A$.  }
\end{defn}

The results of this paper are proved using a version of the
Contracting Mapping Theorem (Theorem \ref{cmt}).  The maps we use
arise from certain vector fields, perhaps defined only locally, on
Riemannian manifolds.  To describe these maps, let $\na$ be the
Levi-Civita connection on a Riemannian manifold $(M,g)$, not assumed
complete. If $X$ is a $C^1$ vector field defined on some open set
$V\subset M$, then at each point $p\in V$ we can view the covariant
derivative $\na X$ as a linear transformation $T_pM\to T_pM$, namely
$v\mapsto \na_v X$.  Call $X$ {\em nondegenerate} on a subset
$U\subset V$ if this endomorphism $(\na X)_p$ is invertible for all
$p\in U$. When referring to bounds on $(\na X)_p^{-1}$ and other
linear transformations, throughout this paper we use the operator
norm: $\|T\|=\sup_{\|v\|=1}\|T(v)\|.$

A $C^1$ vector field $X$ defined on an open set in $M$ and
nondegenerate on a subset $U$ defines a map $\Phi_X:U\to M$
by
\be\label{defphix}
\Phi_X(p)=\exp_p(-(\na X)_p^{-1}X_p) ,
\ee
assuming that $\exp_p(-(\na X)_p^{-1}X_p)$ is defined for all $p\in
U$. (In this paper we use both $X_p$ and $X(p)$ to denote the value of
a vector field $X$ at at a point $p$.) Note that zeroes of $X$ are
fixed-points of $\Phi_X$, and if $\|X\|$ is not too large pointwise
then the converse is true as well.  One of the results of this paper
is the following theorem, a much stronger version of which is proven
in \S 2.

\begin{theorem}\label{thm1}
Let $(M,g)$ be a Riemannian manifold and let $U\subset M$ be open.
Given $\e>0,k_1>0,k_2>0$, let ${\cal X}_{\e,k_1,k_2}(U)$ denote the
set of nondegenerate vector fields $X$ on $U$ satisfying the following
conditions pointwise on $U$: (i) $\|X\|\leq\e$, (ii) $\|(\na
X)^{-1}\|\leq k_1^{-1}$, and (iii) $\|\na\na X\|\leq k_2$. If both $\e
k_1^{-1}$ and $k_2k_1^{-1}$ are sufficiently small, and $X\in {\cal
X}_{\e,k_1,k_2}(U)$, then $\Phi_X:U\to M$ is a contraction, where the
distance function on $U$ is the one determined by the Riemannian
metric $g$ on $M$. If $U$ is a ball $B$ of radius $\rho$ centered at
$p_0$, and if $\rho$ is sufficiently small and $\e,k_1,k_2$ are as
above, then there exists a positive $\e_1\leq
\e$ such that if $\|X(p_0)\|\leq
\e_1$, then $\Phi_X$ preserves $B$ and hence has a unique fixed point
$\Bar{p}$ in $B$; the point $\Bar{p}$ is also the unique zero of $X$
in $B$. For all $p$ in some possibly smaller open ball centered at
$p_0$, the iterates $(\Phi_X)^n(p)$ converge to $\Bar{p}$.
\end{theorem}

\begin{example}\label{example1}
{\em Euclidean space $\bfr^n$.} Since $T_x\bfr^n\iso \bfr^n$
canonically for all $x\in\bfr^n$, a vector field $X$ on $\bfr^n$ can
be naturally identified with a vector-valued function
$G:\bfr^n\to\bfr^n,$ and the Levi-Civita connection is just given by
ordinary directional differentiation: $(\na
X)_x(v)=(DG|_x)(v)=\frac{d}{dt}G(x+tv)|_{t=0}.$ The exponential map is
given simply by $\exp_x(v)=x+v$. Thus
\ben
\Phi_X(x)=x-(DG|_x)^{-1}(G(x)),
\een
which is exactly the Newton's-method map used in the usual
contracting-mapping proof of the Inverse Function Theorem;  cf.
\cite{LS} \S 4.9. 
\end{example}

\ss
Example \ref{example1} illustrates the close relationship between the
iteration in Theorem \ref{thm1} and Newton's method. However, one
gains considerable flexibility by not requiring quite so strict a
relationship as in (\ref{defphix}), looking more generally at maps of
the form $p\mapsto\exp_p(Y_p):=\Psi_Y(p)$ for suitable vector fields
$Y$. Our approach will focus on maps of this more general form,
deducing consequences for maps of the form $\Phi_X$ as a special
case. For the maps $\Psi_Y$, the size restriction on $\|\na X\|$ and
$\|\na \na X\|$ is replaced by the single condition that at each point
the endomorphism $\na Y$ be close to minus the identity. Note that in
this case, $-(\na Y)^{-1}Y$ is close to $Y$, so that the maps $\Phi_Y$
and $\Psi_Y$ are themselves close. Iterative schemes based on maps of
the form $\Psi_Y$ are thus a natural generalization of Newton's
method.  Our most general result for these maps and their associated
algorithms is Theorem
\ref{quant_cor}, a stronger version of Theorem \ref{thm1} in which all
the ``sufficiently smalls'' are quantified for the maps $\Psi_Y$ and
$\Phi_X$. One corollary is the following:

\begin{cor} \label{easycor_intro}
Let $\d\leq\Delta\in\bfr, r_1\in\bfr,$ and
suppose that the sectional curvature $K$ of $M$ satisfies $\d\leq K
\leq \Delta$. There exists a number $D_{\rm crit}$, depending only on
$\d,\Delta,$ and $r_1$, such if $\mu$ is a probability distribution
supported on a set $Q\subset M$ of diameter less than $D_{\rm crit}$,
and the local convexity radius at some point of $Q$ is at least $r_1$,
then the primary center of mass $\Bar{q}$ of $\mu$ exists, and the
Riemannian averaging algorithm converges to $\Bar{q}$ for every
initial point $q\in Q$.
\end{cor}

\ni The definition of  $D_{\rm crit}$ in terms of $\d,\Delta,$ and
$r_1$ is given in \S 4 (see (\ref{defdcrit})); the ``primary center of
mass'' is defined in
\S 3.

\ss
We use the exponential map in defining $\Psi_Y$ because of its
universality, but in specific examples ``exp'' can be replaced by
other maps defined on a neighborhood of the zero-section of the
tangent bundle. This is convenient in the shape-space setting for the
algorithm GPA-S$'$; see \cite{groi2}. However, any continuous map
$F:(U\subset M)\to M$ can always be expressed in the form $\exp\circ
Y$, with $Y$ continuous, provided that for all $p\in U$ the distance
$d(p,F(p))$ is less than the local injectivity radius at $p$ (see
Definition \ref{definj}). Thus if we are interested only in maps that
have any chance of having fixed points, we can always restrict
attention to maps of the form $\Psi_Y$.

This paper is organized as follows. In \S 2 we study the maps $\Psi_Y$
and derive conditions for iterative algorithms based on these maps to
converge.  Before specializing to the Riemannian averaging algorithm,
some discussion of Riemannian centers of mass is needed; this is given
in \S 3, where we also define the vector field $Y$ on which the
averaging algorithm is based.  In general a probability distribution
on a manifold (even one supported on a finite set) can have more than
one center of mass, depending on how ``center of mass'' is defined,
but under certain circumstances one of these is distinguished. In
statistics this is typically done using least-squared-distances
minimization.  However, we offer a more directly geometric way of
singling out a ``primary'' center of mass, using convex hulls. We
digress a bit in Section 3 from the main contracting-mapping theme
because, surprisingly, we have not found any discussion of the
relation of Riemannian centers of mass to convex hulls anywhere in the
center-of-mass literature, although the idea seems very natural.  Our
final statement concerning this relationship, Corollary \ref{primcom},
may be a fact known to workers in the field but it is a stronger
explicit statement than we have seen elsewhere.

In \S 4 we apply the results of \S 2 to obtain a constructive proof of
the existence and uniqueness of the center of mass of a probability
distribution $\mu$ with sufficiently support in a ball of sufficiently
small radius $\rho$ (Corollary \ref{weakkarch}).  Karcher's
existence/uniqueness theorem has a less stringent requirement on
$\rho$, and its uniqueness statement has been strengthened by
W. S. Kendall
\cite{wkendall1}. In view of these results, the most important
feature of the contracting-mapping approach to the center-of-mass
problem is not that it gives existence and uniqueness of the average,
but that it provides a constructive algorithm for finding it (Theorem
\ref{mythm}), along with convergence-rate estimates.  The restriction
on $\rho$ in Theorem \ref{mythm} is almost certainly not sharp.  If
the map on which the algorithm is based has a certain convexity
property that we call ``tethering'', then the upper limit on $\rho$
can be increased considerably. Tethering may occur fairly generally,
but the author has no proof of this. Thus the results in sections 4--6
are stated both without and with the assumption of tethering.

In \S 5 we estimate the convergence rate of algorithms of the form
``iterate $\Psi_Y$'' for general $Y$, and show that the rate is
completely controlled by bounds on $\na Y+I$. In general the
convergence of the sequence $\{p_n=\Psi_Y(p_0)\}$ is exponential; if
$\|\na Y+I\|\leq \e_1$ then $d(p_{n+1},p_n)\leq
d(p_1,p_0)\e_1^{n}$. For maps of the form $\Phi_X$ the convergence is
much faster, obeying the same bounds that one has for Newton's method
in Euclidean space. For the Riemannian averaging algorithm we obtain
something in between: exponential convergence, but with a constant
$\e_1$ that is $O(D^2)$, where $D$ is the diameter of the support of
the distribution being averaged.  We also combine the convergence-rate
result with W. S. Kendall's uniqueness result to obtain a sharpening
of Theorem \ref{mythm} (Theorem \ref{strongerthm2}), establishing
convergence of the algorithm under a weaker requirement on $\rho$.

The statement that $\e_1$ is $O(D^2)$ heuristically---and {\em only}
heuristically---explains the rapid convergence of the GPA algorithms;
it does not fully explain why GPA algorithms converges rapidly in
any applications (or determine in advance whether they will), since
asymptotics do not tell us how small $D$ must be before the leading
asymptotic term decently approximates the actual convergence
rate. However, Theorem
\ref{strongerthm2} can be used to give bounds on $\e_1$ of the form
$\e_1\leq cD^2$ (for all $D$ less than the critical diameter in the
theorem, not just for small $D$), where $c$ is computable from the
geometry of $M$.  In \S \ref{examples} we carry this out and give a
{\em universal} worst-case estimate of the convergence rate when the
curvature of $M$ is non-negative, which is the case in all shape space
and size-and-shape space applications.

In the appendix (\S \ref{appendix}) we prove (or cite proofs of)
certain facts used in \S\S \ref{zvf}--\ref{construct} concerning
Jacobi fields and the distance function.


\setcounter{equation}{0}
\section{Zeroes of Vector Fields}\label{zvf}

Throughout this paper, $M$ denotes a smooth connected manifold
equipped with a Riemannian metric $g$. The induced distance function
on $M\times M$ is denoted $d_M(\cdot,\cdot)$, or simply
$d(\cdot,\cdot)$ when no ambiguity can arise. $M$ is always regarded
as a metric space with this distance function, and the closure of a
subset $U$ in $M$ is denoted $\Bar{U}$.  $B_\rho(p)\subset M$ denotes
the open ball of radius $\rho$ centered at $p$.  If $U\subset M$ is
connected, $d_U$ denotes ``distance within $U$'', the infimum of
lengths of curves in $U$ connecting two given points of $U$.  $TM$
denotes the tangent bundle of $M$, and $\pi:TM\to M$ the canonical
projection. $X$ and $Y$ denote vector fields on $M$ that are at least
$C^2$ and $C^1$ respectively.  If $N_1,N_2$ are manifolds and
$F:N_1\to N_2$ is a smooth map, then for $p\in N_1$, we let
$F_{*p}:T_pN_1\to T_{F(p)}N_2$ denote the derivative of $F$ at
$p$. The identity map of any space is denoted $I$.

The main theorems of this paper are deduced from the following
corollary of the standard Contracting Mapping Theorem (cf. \cite{LS}
Corollary 4.9.2). 

\begin{thm}[Contracting Mapping Theorem] \label{cmt} 
Let $B=B_\rho(p_0)$ be an open ball in a metric space $(A,d)$, with
$(\Bar{B},d)$ complete. Suppose that $B\subset U\subset A$, that
$F:U\to A$ is a contraction with constant $\k$, and that
$d(p_0,F(p_0))<(1-\k)\rho$. Then $F$ preserves $B$  and has a unique
fixed point $\Bar{p}$. Furthermore $\Bar{p}\in B$ and
$\lim_{n\to\infty}F^n(q)=\Bar{p}$ for all $q\in B$.
\end{thm}

As in the Euclidean case (Example \ref{example1}), in the general case
$\Phi_X$ (and more generally $\Psi_Y$) turns out to be a contraction
on sets on which $\|X\|$ (more generally $\|Y\|$) is sufficiently
small.  Our proof of this fact relies on the following simple fact.

\begin{lemma}\label{diffble_contr}
Let $U,M$ be connected Riemannian manifolds and let $\k<1$.
If $F:U\to M$ is a $C^1$ map satisfying
\be\label{constk}
\|F_{*p}\|\leq \k \ \ \mbox{\rm for all\ } p\in
M
\ee
then $F$ is a contraction with constant $\k$.  
\end{lemma}

\pf For any curve $\g$ in $U$ connecting $p$ to $q$, (\ref{constk})
implies $\ell(F\circ\g)\leq \k\ell(\g)$, where $\ell$ denotes
arclength. 
\qedns

We will prove that $\Phi_X$ is a contraction (on suitable sets) by
computing its derivative and applying Lemma \ref{diffble_contr}.  The
map $\Phi_X$ is of the form $\exp\circ Y$, where $Y$ is a vector field
on $M$. Below we express the derivatives of the maps $Y:M\to TM$ and
$\exp:TM\to M$ in terms of the horizontal-vertical splitting of
$T(TM)$ induced by the Levi-Civita connection $\na$. We first review
this splitting (see also \cite{karcher}, Appendix B).

Given a curve $\g$ in $M$ starting at a point $p$ (i.e. a map $\g$
from some interval of the form $(-\e,\e)$ to $M$ with $\g(0)=p$), a
{\em lift} of $\g$ starting at $w\in T_pM$ is a curve $\tg$ with
$\pi\circ\tg=\g$ and $\tg(0)=w$---i.e. a vector field along $\g$ whose
value at $p$ is $w$.  A lift $\tg$ is {\em horizontal} if this vector
field is parallel ($\na_{\g'(t)}\tg\ident 0$). Every curve $\g$ has a
unique horizontal lift starting a a given $w\in T_{\g(0)}M$, and the
vector $\tg'(0)\in T_w(TM)$ depends only on $\g'(0)$. Hence the map
$\g'(0)\mapsto \tg'(0)$ is well-defined and at each $w\in TM$ uniquely
determines a {\em horizontal lift} $\tilde{v}\in T_w(TM)$ of each
$v\in T_{p}M$, where $p=\pi(w).$ The {\em horizontal subspace} of
$T_w(TM)$ is defined to be the subspace $H_w$ consisting of all
horizontal lifts to $w$ of vectors in $T_{p}M$, and
$\pi_{*w}|_{H_w}:H_w\to T_{p}M$ is an isomorphism. The {\em vertical
subspace} $V_w$ of $T_w(TM)$ is the tangent space to the fiber $T_pM$
at $w$. The subspace $V_w$ is canonically isomorphic to $T_pM$
(identifying a vertical vector $\frac{d}{dt}(w+tv)|_{t=0}$ with $v$);
we denote the inverse of this isomorphism by $\i:T_pM\to V_w(TM).$ The
horizontal and vertical subspaces provide a splitting of $T_w(TM)$:
for every $u\in T_w(TM)$, there exist unique vectors $a,b\in T_pM$
such that $u=\tilde{a}+\i(b)$ (specifically $a=\pi_*w$ and
$b=\i^{-1}(w-\tilde{a}));$ we write $\tilde{a}=\hor(u)$ and
$\i(b)=\vert(u)$.

The derivatives we need will be expressed in terms of Jacobi fields
(vector fields $J$ along geodesics $\g$ satisfying the Jacobi equation
\be\label{jaceq}
\na_{\g'}\na_{\g'}J={\rm Riem}(\g',J)\g';
\ee
see \cite{CE} \S 1.4 or \cite{karcher} Appendix A).
Below, for $w,a,b\in TM$ with the same base-point, let $J^w_{(a,b)}$
denote the Jacobi field $J$ along $\g_w$ with $J(0)=a,
(\na_{\g_w'}J)(0)=b$, where $\g_w$ is the unique geodesic with initial
velocity $w$ (i.e. $\g_w(t)=\exp_{\pi(w)}(tw)$).  

Throughout this paper we will be concerned with maps of the form
\be\label{defpsiy}
\Psi_Y=\exp\circ Y :U\to M,
\ee
where $Y$ is a vector field on some domain $U\subset M$. In
(\ref{defpsiy}) we view $Y$ as a map $U\to TM$ and assume that ${\rm
image}(Y)\subset{\rm domain}(\exp)$.  The derivative $(\Psi_Y)_*$ is
given by
\be\label{psiyderiv}
(\Psi_Y)_{*p}v= J^{w}_{(v,0)}(1)+(\exp_p)_{*w}(\i((\na_vY)_p)) \ 
\in T_{\Psi_Y(p)}M
\ee
where $w=Y_p$; the formula above can be deduced from 
\cite{karcher} Appendix B.  If $X$ is a nondegenerate vector field 
and we define $\Phi_X:U\to M$ as in (\ref{defphix}), then for the
vector field $Y=-(\na X)^{-1}X$ we have
\be\label{nay}
\na Y=(\na X)^{-1}\circ(\na\na X)\circ(\na X)^{-1}X-I.
\ee 
Thus as a particular case of (\ref{psiyderiv}) we have
\be\label{derivphix}\label{phixderiv}
(\Phi_X)_{*p}(v)=J^{Y_p}_{(v,0)}(1)-(\exp_p)_{*Y_p}(\i(v))+(\exp_p)_{*Y_p}
(\i((\na X)^{-1}\circ(\na_v\na X)\circ(\na X)^{-1}X)), 
\ee 
where $X,\na X,$ and $\na_v\na X$
are evaluated at $p$.

The term $(\exp_p)_{*Y_p}(\i(v))$ in (\ref{phixderiv}) is itself the
value of a Jacobi field, namely
$J^{w}_{(0,v)}(1)$ where $w=Y_p$. Hence (\ref{derivphix}) can be written as

\be\label{derivphix2}\label{phixderiv2}
(\Phi_X)_{*p}(v)=\hat{J}^{p}_v(1) + (\exp_p)_{*Y_p}(\i(Z_p)) 
\ee 

\ni where $Z_p=(\na X)_p^{-1}(\na_v\na X)_p(\na X)_p^{-1}X_p$ and
where $\hat{J}^p_v$ is the Jacobi field along $\g_{w}$ with the
``antidiagonal'' initial conditions
$\hat{J}(0)=-(\na_{\g'}\hat{J})(0)=v$.  In Euclidean space this
Jacobi field always vanishes at time 1, and $(\exp_p)_*$ is the
identity after appropriate identifications are made as in Example
\ref{example1}, so that (as is well known) $\Phi_X$ is a contraction
if at each point $\|X\|$ is small enough in terms of $\|(\na
X)^{-1}\|$ and $\|\na\na X\|$.  In the general case we can again make
$\|(\exp_p)_*(\i(Z_p)) \|$ arbitrarily small by taking $\|X_p\|$
sufficiently small. Additionally, $\|X_p\|$ small implies $\|Y_p\|$
small, implying that the geodesic $\g_{Y_p}$ is short. For
sufficiently short geodesics, the map $v\mapsto \|\hat{J}^p_v(1)\|$ is
arbitrarily close to the corresponding map on Euclidean space, namely
the zero map.  (We will prove a stronger version of this fact in Lemma
\ref{shortgeods} below.) Hence it is already clear that if
$\sup_p{\|X_p\|}$ is sufficiently small on a set $U$, then
$(\Phi_X)|_U$ will be a contraction.

The essential ingredient in the preceding argument is that $\Phi_X$ is
a map of the form $\Psi_Y=\exp\circ Y$ for some vector field $Y$ whose
covariant derivative is close to minus the identity (pointwise)
whenever $\|Y\|$ is small enough. (The prototypical example is the
radial vector field $-\sum_ix^i\frac{\partial}{\partial x^i}$ on
$\bfr^n$, whose covariant derivative is identically $-I$.)  In
computational situations it may be costly to invert $\na X$, so we
will analyze the more general maps $\Psi_Y$, and deduce results for
maps of the form $\Phi_X$ as a special case.

For some applications (e.g. those in \cite{groi2}), it is useful to
know the explicit dependence of our eventual contraction constants on
background geometric parameters, so we  keep track of this
dependence carefully---leading unavoidably to longer formulas than if
we were only aiming at qualitative results. Certain special functions
will appear, all of which are related to the analytic (entire)
functions $\ch,\sh$ defined by
\be\label{chsh}
\ch(z)=\sum_{n=0}^\infty \frac{z^n}{(2n)!}, \ \ \ 
\sh(z)=\sum_{n=0}^\infty \frac{z^n}{(2n+1)!}.
\ee
Since the definitions and properties of the relevant functions are
scattered through the text, for reference Table 1
lists the functions and the properties used.

\begin{table}[hp] \label{table1}
\begin{center}{Table of Special Functions}\\ 
\ss
\begin{tabular}{|c|c|c|}
function & defining formula & properties used\\
\hline
$\ch(z),\ z\in\bfc$& $\sum_{n=0}^\infty z^n/(2n)!$  & \\
\hline
$\sh(z),\ z\in\bfc$ &$ \sum_{n=0}^\infty z^n/(2n+1)!$& \\
\hline
$\phi_-(x),\ x\in[0,\infty)$ & 
$\begin{array}{ll}
\ch(x^2)-\sh(x^2)\\ =\cosh x- x^{-1}\sinh x 
\end{array}$& mono. $\uparrow$, $\phi_-(x)\geq 0$\\
\hline
$\phi_+(x),\ x\in[0,3\pi/4)$ & 
$\begin{array}{ll}
\sh(-x^2)-\ch(-x^2)\\ = x^{-1}\sin x-\cos x
\end{array}$
&
mono. $\uparrow$, $\phi_+(x)\geq 0$\\
\hline
 $C_1(\l,r), \l\in\bfr, r\geq 0$
&
$
\begin{array}{cl}  1, & {\rm if}\ \l\geq 0\\
\frac{\sinh(|\l|^{1/2}r)}{|\l|^{1/2}r},  {\rm if}\ & \l<0
\end{array}$
& \begin{tabular}{l}
mono. $\uparrow$ in each variable, \\
 $C_1(\l,r)\geq 1$
\end{tabular}\\
\hline
\begin{tabular}{l}
$h(\l,r),\ \l\in\bfr, $
\\
\mbox{\hspace{.5in}}
$r\in 
\left\{ 
\begin{array}{l} 
[0,\pi) \ {\rm if}\ \l>0,  \\  
{[0,\infty) } \ {\rm if}\ \l\leq 0 
\end{array}
\right.
$
\end{tabular}
& 
$\ch(-\l r^2)/\sh(-\l r^2)$
&
\begin{tabular}{l}
$h(0,r)=h(\l,0)=1$, \\
$h(\l,r)>0$ if $\l\leq 0$,\\
or if $\l>0$ and $\l^{1/2}r<\pi/2$
\end{tabular}
\\
\hline
$h_-(x)=h(-1,x), \ x\in [0,\infty)$&   $x\coth x$ & 
\begin{tabular}{l}mono. $\uparrow$, \\$h_-(x)\geq h_-(0)=1$
\end{tabular}\\
\hline
$h_0(x)=h(0,x), \ x\in [0,\infty)$ & $1$ & \\
\hline
$h_+(x)=h(1,x), \   x\in [0,\pi)$ & $x\cot x $ & 
\begin{tabular}{l} mono. $\downarrow$,\\ $h_+(x)\leq h_+(0)=1$
\end{tabular}\\
\hline
$\psi(\l,r),$ same domain as $h$
&
${\rm sign}(\l)(1-h(\l,r))$ & 
\begin{tabular}{l} 
$\psi(\l,r)\geq 0$, mono. $\uparrow$ \\
in $|\l|$ and $r$, \\ convex in each variable
\end{tabular}
\\
\hline
$\begin{array}{l}
\psimax(\d,\Delta,r), \ \d\leq \Delta\in\bfr, \\
\mbox{\hspace{.75in}} r\in[0,\infty)
\end{array}$
& $\max(\psi(\Delta,r),\psi(\d,r))$ & 
\begin{tabular}{l}
mono. $\uparrow$ in $\Delta$ and $r$, \\
mono. $\downarrow$ in $\d$, 
\\
convex in each variable,\\
$\psimax(\d,\Delta,0)=0$
\end{tabular}\\
\hline

\end{tabular}
\end{center}
\caption{\footnotesize
In this table and throughout this paper our convention for functions
that are given for $x\neq 0$ by formulas such as ``$x^{-1}\sin x)$'' are
extended to $x=0$ by continuity.  When monotonicity or convexity of a
multivariable function is stated with respect to one variable, the
other variables are assumed fixed.  }
\end{table}

\bs
To estimate $\|(\Psi_Y)_*\|$, we rewrite (\ref{psiyderiv}) as
\be\label{psiyderiv2} (\Psi_Y)_{*p}v=
\hat{J}^p_v(1)+(\exp_p)_{*Y_p}(\i((\na Y|_p+I)v)).
\ee

\ni We will first analyze the Jacobi fields $\hat{J}^p_v$. 

\ssn{\bf Notation.}
For any subset $U\subset M$, let $\Delta(U)$ and $\d(U)$ denote,
respectively, the supremum and the infimum of the sectional curvatures
of $(U,g|_U)$; let $\aku(U)=\max(|\Delta(U)|,|\d(U)|)$.  For a curve
$\g$ we simply write $\Delta(\g)$ for $\Delta({\rm Im}(\g))$, etc.
Then we have the following proposition. The inequality (\ref{jbound})
below can be derived from Karcher's elegant (and more general)
Jacobi-field bounds; see \cite{karcher} pp. 534-535, 539. However, for
the special case (\ref{jbound}), we give a short, direct proof in the
Appendix (\S \ref{app2}).  In the second part of \S \ref{app2} we show
how the proof leads directly to (\ref{jbound_symplus}).

\begin{prop} \label{shortgeods}
Let $p\in M$, let $\g:[0,1]\to M$ be a geodesic of length $r$ starting
at $p$, and for each $v\in T_pM$ let $\hat{J}_v$ be the Jacobi field
along $\g$ with the ``antidiagonal'' initial conditions
$(\hat{J}_v(0), (\na_{\g'}\hat{J}_v)(0)) = (v,-v).$ Let $v^\perp$
denote the component of $v$ perpendicular to $\g'(0)$.  Then

\be \label{jbound}
\| \hat{J}_v(1) \| \leq \phi_-(r\aku(\g)^{1/2})\|v^\perp\|
\ee
where
\be\label{defphim}
\phi_-(x)=\cosh(x)-\frac{\sinh x}{x}. 
\ee
If $M$ is a locally symmetric space of nonnegative curvature, and
$\Delta(\g)^{1/2}r<3\pi/4$, this bound can be sharpened to
\be \label{jbound_symplus}
\| \hat{J}_v(1) \| \leq \phi_+(r\Delta(\g)^{1/2})\|v^\perp\|
\ee
where
\be\label{defphip}
\phi_+(x)=\frac{\sin x}{x}-\cos x.
\ee
\end{prop}

The ``$3\pi/4$'' in the locally-symmetric case can be increased to
approximately $.87\pi$ (see the discussion of (\ref{replphi}) in \S
\ref{app2})) but any instances in which $r\Delta(\g)^{1/2}>\pi/2$ are
irrelevant for all uses in this paper.

\ms
Turning our attention to the second term in (\ref{psiyderiv2}), we
have 
\ben
\|(\exp_p)_{*Y_p}(\i(\na Y|_p+I)v)\| \leq 
\|(\exp_p)_{*Y_p}\| \ \|(\na Y|_p+I)\| \ \|v\|.
\een

\ni We recall the following terminology.

\begin{defn} \label{definj}{\rm
The {\em local injectivity radius} at $p\in  M$ is
$\rinj(p):=\sup\{\rho \mid \exp_p: (B_\rho(0)\subset T_pM)\to M$ is
defined and is a diffeomorphism onto its image\}; $\rinj(\cdot)$ is a
positive continuous function on $M$.  For any subset $U\subset M$, we
define $\rinj(U)=\inf_{p\in U}\{\rinj(p)\}$. When $U=M$ this infimum
is called the {\em injectivity radius of $(M,g)$}.}
\end{defn}

\begin{defn}\label{defcvx}{\rm
 A subset $U\subset M$ is {\em convex} (respectively, {\em strongly
convex}) if for all $p,q\in U$ (resp., for all $p\in U, q\in \Bar{U}$)
there is a unique minimal geodesic segment $\g$ in $M$ from $p$ to
$q$, and $\g-\{q\}$ lies entirely in $U$\footnote{In the differential
geometry literature there is little consistency in the meanings
attached to the terms ``convex set'' and ``strongly convex set''.
There is quite an array of criteria one can imagine demanding of a
convex set; see Definition \ref{vis} for a few of these.}.  For each
$p\in M$ we define the {\em local convexity radius} $\rcvx(p):=\sup\{
\rho\leq
\rinj(p)
\mid B_\rho(p)\ \mbox{\rm is convex}\}$; for $U\subset
M$ we let $r_{\rm cvx}(U)=\inf_{p\in U}\{r_{\rm cvx}(p)\}$.  Like the
local injectivity radius, the local convexity radius of a point (or of
a closed set) is always positive (\cite{helgason} Lemma I.6.4). } 
\end{defn}

Convexity is relevant because we want to apply Lemma
\ref{diffble_contr} and Theorem \ref{cmt} to the case $U\subset M$.
The lemma only gives us a contraction from the metric space $(U,d_U)$
to $(M,d_M)$. However if $U$ is convex then $d_U=d_M$ so that Theorem
\ref{cmt} applies.

For $w\in T_pM$ with $\|w\|<\rinj(p)$ the norm of $(\exp_p)_{*w}$ can
be bounded in terms of curvature and $\|w\|$:

\be
\|\exp_{p*w}\| \leq C_1(\d(\g),\|w\|)
\ee 
where $\g$ is the geodesic from $p$ with $\g'(0)=1$ and where
\be\label{defc2}
 C_1(\l,r)=\left\{
\begin{array}{cl}  1, & \l\geq 0\\
\frac{\sinh(|\l|^{1/2}r)}{|\l|^{1/2}r}, & \l<0
\end{array}\right. 
\ee
(see \cite{karcher} estimate C1).
\ni 
Thus if the image of $\g$ lies in a set $U$, and $\|Y_p\|<\rinj(p)$, then
\be\label{2ndterm} 
\|(\exp_p)_{*Y_p}(\i(\na Y|_p+I)v)\| \leq
C_1(\d(U), \|Y_p\|)\ \|(\na Y +I)_p\| \ \|v\|. 
\ee 
Assembling the pieces above, we have the following corollary.

\begin{cor}\label{qual_cor}
Let $(M,g)$ be a Riemannian manifold, $\rho>0$, $p\in M$, and
$B=B_\rho(p)$. Assume that $\rho\leq\rcvx(B)$ and that
$\aku(B)<\infty$.

(a) There exists $\e>0$, depending only on $\rcvx(B)$ and the
sectional curvature of $(B,g)$, such that if $Y$ is a
vector field defined on $B$, with $\|Y\|\leq \e$
and $\|\na Y + I\| \leq \e$ pointwise on $B$, then
$Y$ has a unique zero in $B$, namely $\lim_{n\to\infty} (\Psi_Y)^n(q)$
for any $q\in B$.

(b) Let $k_1,k_2 >0$. There
exists $\e>0$, depending only on $k_1, k_2,\rcvx(B)$, and the
sectional curvature of $(B,g)$, such if $X$ is a vector field $X$
satisfying $\|(\na X)^{-1}\| < k_1$, $\| \na\na X\| \leq k_2$, and 
$\|X\|\leq \e$ pointwise on $B$, then $X$ has a unique zero in $B$, namely
$\lim_{n\to\infty} (\Phi_X)^n(q)$ for any $q\in B$.
\end{cor}

\begin{remark} \label{inc}{\rm We intentionally avoid assuming that
that $(M,g)$ is complete or has positive injectivity radius. In the
application to the set of smooth points of the shape space
$\Sigma^k_n$, if $n\geq 3$ then $(M,g)$ is a dense open subset of a
non-smooth real algebraic variety (cf. \cite{carne}), hence neither
complete nor of positive injectivity radius. However, any closed
subset of $M$ with positive injectivity radius will be complete.  In
particular this applies to the closures of all the balls considered in
this paper.}
\end{remark}

\ss
Corollary \ref{qual_cor} follows immediately from the following
more quantitative version.

\begin{thm}\label{quant_cor}
Let $U\subset M$ be connected and let $\aku=\aku(U), \d=\d(U)$. Define
the functions $\phi_-(\cdot)$ and $C_1(\cdot,\cdot)$ by {\rm
(\ref{defphim})} and {\rm (\ref{defc2})}.  Assume either of the
following sets of hypotheses:

\ss
{\em Case 1.} 
$Y$ is a vector field defined on $U$ and at each
point of $U$ we have $\|Y\|\leq \e_0<\rinj(U)$ and $\|\na Y + I\| \leq
\e_1.$ Define $\Psi_Y=\exp\circ Y$ as in (\ref{defpsiy}).

\ss
{\em Case 2.}
$k_1,k_2>0$, $X$ is a vector field defined and uniformly nondegenerate on
$U$, and at each point of $U$ we have
$\|(\na X)^{-1}\|\leq k_1^{-1}$, $\|\na\na X\|\leq k_2$, and $\|X\|\leq
\e<k_1\rinj(U).$ Define $\Phi_X=\exp\circ (-(\na X)^{-1}\circ X)$ as
in (\ref{defphix}). 

\ssn Then:

\ss
(a) For all $p\in U$, in Case 1 we have
\be\label{psibound}\label{defkpsiy}
\|(\Psi_Y)_{*p}\| \leq \k(\Psi_Y):=
\phi_-(\aku^{1/2}\e_0) + C_1(\d,\e_0)\e_1,
\ee
while in Case 2
\bearray\label{phibound}\label{defkphix}
\| (\Phi_X)_{*p} \| &\leq& \k(\Phi_X):=
\phi_-(\aku^{1/2}\e k_1^{-1}) + C_1(\d,\e k_1^{-1})k_2k_1^{-2}\e.
\eearray

\ssn In Case 1, let $F=\Psi_Y$; in Case 2 let $F=\Phi_X$.
In Case 1 (respectively Case 2) if $\e_0,\e_1$ are small enough
(resp., $\e$ is small enough) that $\k(F)<1,$ then $F:(U,d_U)\to
(M,d_M)$ is a contraction with constant $\k(F)$, and therefore
has at most one fixed point in $U$.  If $U$ contains an open ball
$B=B_\rho(p_0)$ on whose closure the distance functions $d_U,d_M$
coincide (a condition satisied by every subset of $\Bar{U}$ if $U$ is
convex), and if
\be\label{kplt1}
\|Y(p_0)\|<(1-\k(F|_B))\rho  \mbox{\hspace{.5in} (in Case 1)},
\ee
or 
\be\label{kpslt1}
\|X(p_0)\|<(1-\k(F|_B))k_1\rho  \mbox{\hspace{.5in} (in Case 2)},
\ee
then $F:U\to M$ has a unique fixed point, and this fixed point lies in
$B$. Equivalently, the vector field $Y$ in Case 1, or $X$ in Case 2,
has a unique zero in $U$, and this zero lies in $B$.  Assuming
(\ref{kplt1}) or (\ref{kpslt1}) as appropriate, $F$ preserves $B$, and
the fixed point is $\lim_{n\to\infty}F^n(q)$ for any $q\in B$.

\ss
(b) If $M$ is a locally symmetric space of non-negative curvature
bounded above by $\Delta$,  then in {\rm (\ref{psibound})} and {\rm
(\ref{phibound})}, we can replace the right-hand sides by the smaller
bounds
\bearray
\k_{{\rm sym}+}(\Psi_Y) &:=& 
\phi_+(\Delta^{1/2}\e_0)+\e_1
\\
{\rm and} \ \ \ \
\k_{{\rm sym}+}(\Phi_X) &:=& \phi_+(\Delta^{1/2}\e k_1^{-1})
+k_2k_1^{-2}\e
\eearray
respectively, provided $\Delta^{1/2}\e_0<3\pi/4$ in the first case and
$\Delta^{1/2}\e k_1^{-1}\leq 3\pi/4$ in the second. 

\end{thm}

\pf (a) Case 1. The bound (\ref{psibound}) follows from 
Proposition \ref{shortgeods}, and (\ref{2ndterm}).  If $\k(\Psi_Y)<1$
Lemma \ref{diffble_contr} implies that $\Psi_Y:(U,d_U)\to(M,d_M)$ is a
contraction with constant $\k$. To use the fixed-point theorem we need
a contraction with respect to a single distance
function. However, the assumption that $d_U=d_M$ on $B$ implies that
the restriction of $\Psi_Y$ to $B$ is a $\k$-contraction from
$(\Bar{B},d_M)\to (M ,d_M)$.  As noted in Remark \ref{inc}, the metric
space $(\Bar{B},d_M)$ is complete. Hence the result follows from Theorem
\ref{cmt} with $U=\Bar{B}$.

\ss Case 2. Letting $Y=-(\na X)^{-1}X$, for $p\in U$ we have $\|Y_p\|\leq
\e k_1^{-1}<\rinj(U)$, and from (\ref{nay}) we have $\|(\na Y +I)_p\|\leq
k_1^{-2}k_2.$ Hence Case 2 follows from Case 1.

\ss (b) This follows from (\ref{jbound_symplus})
and the proof of (a).
\qed

\begin{remark}{\rm
In the bound (\ref{phibound}), as either $\e\to 0$ or $\aku\to 0,$ we
have $\phi_-(\aku^{1/2}\e k_1^{-1})\to 0$ and $C_1(\d,\e k_1^{-1})\to
1$. Hence, as one would hope, for small $\e$ and for small $\aku$ the
bound 
(\ref{phibound}) is asymptotic to $k_2k_1^{-2}\e$, the well-known
bound for the Euclidean case (see the discussion following
(\ref{phixderiv2})).}
\end{remark}

\begin{remark}{\rm Theorem \ref{thm1} follows immediately from Case 2
of Theorem
\ref{quant_cor}(a).}
\end{remark}


\setcounter{equation}{0}
\section{Averaging Points in a Riemannian Manifold}\label{3_com}

In its most elementary form, averaging is something that one does to a
finite list of elements in a vector space.  The average of a list
$\{w_1,\dots, w_m\}$ in a vector space $V$ can be uniquely
characterized as that vector $\Bar{w}\in V$ for which
\be\label{property1}
\sum_{i=1}^m (w_i-\Bar{w})=0.
\ee
The ``balancing property'' (\ref{property1}) motivates the alternative
term for the average, {\em center of mass}.  If $V$ is given {\em any}
inner product then, using the the inner product to define a norm, the
average above can also be uniquely characterized as
\ben
\Bar{w}= \mbox{\rm that vector $v$ which
minimizes\ } \sum_{i=1}^m\|w_i-v\|^2
\een
(the ``least-squares property'').

Unlike the balancing property, which requires a linear structure on
$V$, the least-squares property makes sense if $V$ is replaced by any
metric space.  A {\em Fr\'{e}chet mean} of a finite subset of a metric
space $(A,d)$ is an element $a\in A$ at which the function $ p\mapsto
\sum_{q\in Q} d(q,p)^2$ attains an absolute minimum. In general a
Fr\'{e}chet mean need not exist or be unique, but when it exists
uniquely it is not unreasonable to call it the average of $Q$.

Modulo existence and uniqueness, Fr\'{e}chet means give a way to
extend the notion of ``average'' to finite lists of points in a
Riemannian manifold, or more generally probability distributions on
such a manifold. This idea of the {\em Riemannian center of mass}
dates back at least as far as E. Cartan \cite{cartan1} in the case of
simply connected manifolds of nonpositive curvature; in this setting
the Fr\'{e}chet mean of any probability distribution exists uniquely.
However, the arbitrary-curvature case seems not to have been studied
systematically until the 1970's, when it was investigated principally
by Karcher and Grove (\cite{karcher,GK,grove}; see also \cite{hild}
\S\S 4--5).

Unlike in Euclidean space, on a general Riemannian manifold it is
clear that some restriction on the set of points to be averaged is
necessary; for example there is no reasonable way uniquely to define
the average of antipodal points on a sphere.  Averaging can be done
sensibly only on sets satisfying some suitable convexity condition (of
which there are several). One notion of convexity was given in
Definition \ref{defcvx}; some other relevant notions are given
below. The reader is warned that different authors attach different
names to these notions.

\begin{defn}\label{vis}
 {\rm Let $U\subset M$. We call $U$ 
\begin{itemize}
\item {\em self-visible} if any two
points of $U$ can be joined by at least one geodesic, not necessarily
minimal, lying in $U$;
\item {\em simple} if for any two points
in $U$ there is exactly one connecting geodesic lying in $U$;
\item {\em
solipsistically convex} 
if for any two points $p,q\in U$ there exists a connecting geodesic in
$U$ whose length is minimal among all connecting arcs lying in $U$
(hence of length $d_U(p,q)$).  
\end{itemize}
A function $f$ defined on a self-visible
set $U$ is called {\em (strictly) convex on $U$} if its restriction to
every geodesic in $U$ is a (strictly) convex function of the arclength
parameter.}
\end{defn}

If $f$ is $C^2$ then a sufficient condition for $f$ to be convex on
$U$ is that its covariant Hessian be positive-semidefinite on $U$;
strict positivity implies strict convexity.

\begin{defn}{\rm (cf. \cite{hild} p. 3) 
An open ball $B=B_\rho(p)$ is a {\em regular geodesic ball}
if (i) $\rho<\rinj(p)$, and (ii)
\be\label{regball}
\rho\cdot\max(0,\Delta(B))^{1/2}<\pi/2 .
\ee
For $p\in M$ define the {\em regularity radius} 
\ben
\rreg(p):=\sup\{\rho\mid B_\rho(p)\ \mbox{is a regular geodesic
ball}\} 
\een
and the
{\em regular convexity radius}
\ben
\rregcvx(p)=\min(\rreg(p), \rcvx(p)).
\een
}
\end{defn}

For regular geodesic balls one has the following theorem of Jost
\cite{jost}; see \cite{hild} Theorem 5.3 and \cite{wkendall1} 
Theorem 1.7.

\begin{theorem}\label{jostthm}
Let $B$ be a regular geodesic ball in a complete Riemannian
manifold.  Then $\Bar{B}$ is simple and solipsistically convex, and
geodesics in $B$ contain no pairs of conjugate points. 
\end{theorem}

Completeness of the ambient manifold is not essential in Theorem
\ref{jostthm}; if $B=B_\rho(p)$, it suffices that $\exp_p$ be defined
on the closed ball of radius $\rho$ centered at $0\in T_pM$.  The
example of an open ball of radius $\pi$ in the unit circle shows that a
regular geodesic ball need not be convex. More generally Theorem
\ref{jostthm} implies that regular geodesic ball $B\subset M$ is
convex if and only if the distance functions $d_B$ and $d_M$ coincide
on $B$.

There are various relations among $\rinj,\rcvx,$ and $\rreg$; we
mention only a few.  By definition, $\rinj(p)$ is the largest of the
three radii at $p$.  If $M$ is complete and has constant positive
curvature, then Bonnet's Theorem (\cite{CE} Theorem 1.26(2)) implies
that $\rcvx(p)\leq\rreg(p)$.  But in general, a geodesic ball can be
convex but not regular (see \cite{groi2} for an example), or, as the
circle example shows, regular but not convex.

\ssn{\bf Notation.} If $p,q\in M$ can be joined by a unique minimal
geodesic, we denote by $\exp_p^{-1}(q)$ the unique pre-image of $q$
(under $\exp_p$) of smallest norm.  

\ss
Now let $Q$ be an arbitrary subset of a convex set $U\subset M$, and
let $\mu$ be a probability measure on $Q$. 
For each
$p\in U$ define
\bearray
Y_Q(p) &=& \int_Q \exp_p^{-1}(q)\ d\mu(q) \ \in T_pM, \label{defyq}\\
f_Q(p) &=& \frac{1}{2}\int_Q d(p,q)^2\ d\mu(q);\label{deffq}
\eearray

\ni More properly these objects should be subscripted with the pair
$(Q,\mu)$. However, in most of our results $\mu$ enters primarily
through the geometry of $Q$ rather than in the behavior of $\mu$ on
$Q$.  To emphasize this we will stick to the imperfect notation above.

\begin{defn}\label{defcom}{\rm Let $U\subset M$ be convex.
(1) Let $Q\subset U,$ let $\mu $ be a probability measure on $Q$, and
define a vector field $Y_Q$ by (\ref{defyq}).  If $Y_Q(p)=0$ at a
unique point $p\in U$, we call $p$ the {\em (Riemannian) center of
mass} of $(Q,\mu)$, {\em relative to} $U$.  (2) Let
$\tq=\{q_,\dots,q_m\}$ be a finite list of points in $U$, let $Q$ be
the set of distinct elements of $\tq$, let $\mu$ be the normalized
counting-measure on $Q$, and define $Y_Q$ as above. If $Y_Q(p)=0$ at a
unique point $p\in U$, we call $p$ the {\em Riemannian average} of the
list $\tq$, {\em relative to} $U$.

We call a point \underline{a} center of mass of $(Q,\mu)$
(respectively, \underline{a} Riemannian average of the list $\tq$) if
it is the center of mass of (resp., Riemannian average) relative to
some convex superset.  }
\end{defn}

For a finite list $\tq$, the definition of Riemannian average relative
to $U$ is simply the zero (assumed unique in $U$) of the vector field
$Y=Y_\tq$ on $U$ defined by $Y(p) =\frac{1}{m}\sum_i
\exp_p^{-1}(q_i) \in T_pM.$ Since $\sum_i (\exp_p^{-1}(q_i)-Y_p)=0$,
heuristically, $Y(p)$ represents ``balanced'' average of the points
$q_i$ as seen from $p$.  Alternatively, we can define $f_\tq:U\to\bfr,
f_\tq(p)=\frac{1}{2m}\sum_{i=1}^m d(q_i,p)^2,$ and assume that $f_\tq$
is minimized uniquely at $\bar{q}\in U$. The Gauss Lemma (\cite{CE}
p. 8]) implies that $\grad(d(q,\cdot))|_p=-\exp_p^{-1}(q)$, so
$\grad(f_\tq)=-Y_\tq,$ implying that $Y_\tq$ has its zero at
$\bar{q}$.  Hence Definition \ref{defcom} extends both the
``balancing'' and ``least-squares'' properties of the Euclidean
average.

\begin{remark} \label{bary}{\rm 
Definition \ref{defcom} generalizes easily to a
solipsistically convex or simple set $U$. In this case denote by
$\exp_p^{-1,U}(q)$ that pre-image $v$ of $q$ (under $\exp_p$) of
smallest norm for which $\exp_p(tv)\in U, 0\leq t\leq 1$.  In
(\ref{defyq}) we can replace $\exp_p^{-1}(q)$ by $\exp_p^{-1,U}(q)$,
and $d(p,q)$ by $\|\exp_p^{-1,U}(q)\|$ (in the solipsistically convex
case this is just $d_U(p,q)$). With these replacements it is
still true that $\grad(f_Q)=-Y_Q$, but the interpretation of $Y_Q(p)$ as
an average of points as seen from $p$ is less compelling.  }
\end{remark}

We will refine  Definition \ref{defcom} later for a case in which
one center of mass is singled out, allowing us to dispense with the
awkward ``relative to $U$'' (Definition \ref{defcom2}).

Following \cite{wkendall1,Le1,Le2}, for example, we will call any
relative minimum of $f_Q$ a {\em Karcher mean}. Thus a Fr\'{e}chet
mean is necessarily a Karcher mean, but, absent extra hypotheses, not
vice-versa.  A center of mass of $(Q,\mu)$ under Definition
\ref{defcom} is simply a Karcher mean
that lies inside some convex superset of $Q$.

Karcher proves a somewhat more general version of the following
theorem (\cite{karcher} Theorem 1.2, Definition 1.3, and Theorem 1.5).

\begin{theorem}[Karcher]\label{karchthm}
Let $(M,g)$ be a Riemannian manifold. Assume that $Q\subset B\subset
M$, where $B=B_\rho(p_0)$ is a strongly convex ball. Let
$\Delta=\Delta(B)$ be the supremum of the sectional curvatures in
$B$. Then, with $f_Q$ and $Y_Q$ defined as above, 

(a) $\grad(f_Q)=-Y_Q$.

(b) The function $f_Q$ achieves a minimum value on
$B$, and hence $Y_Q$ has a zero in $B$.

(c) If
$\rho\cdot\max(0,\Delta(B))^{1/2} <\pi/4$, then the minimum of $f_Q$
on $B$ is achieved at a unique point $\Bar{q}$, and for any point
$p\in B$ we have
\be\label{karchdist}
d(p,\Bar{q})\leq \|Y_Q(p)\|\cdot 
\left\{\begin{array}{ll} 1/h(\Delta,2\rho) & {\rm if} \ \Delta>0\\
1 & {\rm if} \ \Delta\leq 0\end{array}\right.
\ee
where, for $\Delta>0$, $ h(\Delta,x) =\Delta^{1/2}x
\cot(\Delta^{1/2}x).$

\end{theorem}

In \cite{karcher}, Karcher defines the center of mass to be the
location of the minimum of $f_Q$ on $\overline{B_\rho}$.  However, his
proof of existence and uniqueness of the minimum also implies
uniqueness of the zero of $Y_Q$, so under the hypotheses of Theorem
\ref{karchthm} this definition coincides with ours; indeed, the
geometric Definition \ref{defcom} is the one used in
\cite{GK}.

Note that the ball $B$ in Theorem \ref{karchthm} is both geodesically
convex and regular. If $\rho<\frac{1}{2}\rreg(p_0)$ then the
requirement on $\rho$ in (c) is automatically satisfied; hence the
upper limit on the radius of the balls for which part (c) is
applicable is at least $\min(\frac{1}{2}\rreg(p_0),\rcvx(p_0))$ but no
greater than $\rregcvx(p_0)$.  In \cite{wkendall1} Theorem 7.3,
W. S. Kendall strengthened the uniqueness assertion\footnote{Kendall's
proof does not yield existence of Karcher means as we have defined
them.  It is clear from the context and the proof that the existence
asserted in the theorem as stated in \cite{wkendall1}
is the existence of a ``solipsistic Karcher mean'', in which the distance
function $d=d_M$ in (\ref{deffq}) is replaced by $d_B$.  The existence
argument requires $\grad(f_Q)$ to be outward-pointing on the boundary
of the ball, which is guaranteed only under the solipsistic
interpretation of $f_Q$.} in Theorem \ref{karchthm}(c):

\begin{theorem}\label{kendthm}{\rm \bf (W. S. Kendall)}
A mass distribution supported in a regular geodesic ball $B$ has at most
one Karcher mean in $\Bar{B}$.
\end{theorem}

In other words, as far as the uniqueness statement is concerned, as
long as we assume $\rho<\rinj(p_0)$ Karcher's $\pi/4$ can be
replaced with $\pi/2$, and the ball $B_\rho(p_0)$ need not be assumed
convex.

In general, Karcher means are not unique in the large, cf.
\cite{CK,wkendall2}. For example, given a set $Q$ of two
equally-weighted points in the unit circle $S^1$, the midpoints of
each of the two arc joining the points is a Karcher mean.  The
statistically-natural absolute minimization of $f_Q$ of course
distinguishes one of these midpoints as the preferred one.  However,
we suggest an alternative, purely geometric way of distinguishing one
of the Karcher means from the rest: just as in Euclidean space, the
center of mass of a distribution $\mu$ should be in the convex hull,
suitably defined, of its support---the average of a set $Q$ should be
not only near $Q$, but ``within'' $Q$.  In the $S^1$ example above,
unless the two points are antipodal---in which case the convex hull is
not defined---only one of the two midpoints meets this criterion.
Thus in this example the convex-hull and global-minimization criteria
coincide, but the author does not know to what extent these criteria
overlap in general.

The definition of ``convex hull'' varies in the literature. The notion
best tailored to our needs is that of the {\em o-hull} defined below.

\begin{defn}\label{defhull}{\rm
Call a set $Q\subset M$ {\em hulled} if it is contained in some convex
set, and {\em o-hulled} if it is contained in some {\em open} strongly
convex set. If $Q$ is hulled (resp. o-hulled), define the {\em convex
hull} of $Q$ (respectively, the {\em convex o-hull} of $Q$), written
$\hull(Q)$ (resp., $\ohull(Q)$) to be the intersection of all convex
sets (resp. open strongly convex sets) containing $Q$. We will usually
refer to these objects just as hulls and
o-hulls.

Note that if a set is hulled, then the minimal geodesic between any
two of its points exists and is unique.  

}
\end{defn}

Obviously hulls and o-hulls, when they exist, are convex sets, and
$\hull(Q)\subset\ohull(Q)$. The o-hull may fail to exist even
when the hull exists (example in $S^1$:
a semicircle closed at one endpoint and open at the other). However 
in $\bfr^n$, at least, the differences between hull and o-hull are
minor: one always has 
\be\label{hoh}
\hull(Q)\subset\ohull(Q) \subset\Bar{\hull(Q)}
\ee
(both inclusions can be strict; see \cite{groi1}). Conceivably
(\ref{hoh}) holds generally for o-hulled sets in Riemannian manifolds
provided ${\hull(Q)}$ has compact closure.

All sets $Q$ of interest in this paper are contained in a convex open
ball and so are o-hulled.  As noted above, we will use o-hulls to
distinguish one particular center of mass.  Neither Karcher's theorem
nor Kendall's generalization, as stated, immediately eliminates the
unsettling possibility that $Q$ could be contained in two different
convex regular geodesic balls, and that $Y_Q$ could have two zeroes
(each of which could even be an absolute minimum of $f_Q$), each
contained in one ball but not the other.  However, the proofs in
\cite{karcher} and \cite{wkendall1} imply more than is explicitly 
stated in either paper, and a minor extension of an ingredient of
these proofs shows that this unwanted phenomenon cannot
happen\footnote{\cite{Le1} uses a different partial solution to this
problem: if in Karcher's theorem it is additionally assumed that
$2\rho<\rinj(p_0)$ and hypothesis (c) is satisfied with $\Delta(B)$
replaced by $\Delta(B_{2\rho}(p_0))$---then Kendall's theorem implies
that the Karcher mean of $(Q,\mu)$ in $B_\rho$ is the unique
Fr\'{e}chet mean of $(Q,\mu)$. Our alternative approach does not
require this extra hypothesis in order to single out a ``best''
Karcher mean, but our geometric definition of ``best'' differs from
the statistical definition.}. We give this extension in Lemma
\ref{unimin} and Corollary \ref{onecom}.  The corollary leads us to
the convex-hull criterion in Definition
\ref{defcom2} below.  

While $\ohull(Q)$ is the {\em smallest} set we can construct naturally
from the family of open strongly convex supersets of $Q$, the {\em
largest} set we can construct from this family also has relevance:

\begin{defn}\label{defstar}{\rm
For any o-hulled set $Q\subset M$, define 
$\Star(Q)$ to be the union of all open strongly convex
supersets of $Q$. Analogously, define $\regstar(Q)$ to be the union of
all regular geodesic balls containing $Q$. Note that $\Star(Q)$
depends only on $\ohull(Q)$. }
\end{defn}

Given an open set $U\subset M$ and a boundary point $p\in \partial U$,
call a tangent vector $v\in T_p M$ {\em outward-pointing} for $U$ if
$v\neq 0$ and if for some $C^1$ curve in $M$ with $\g'(0)=-v$ we have
$\g((0,\e))\subset U$ for some $\e>0$. 

For reference, we record the following obvious facts (proof left to
the reader).

\begin{lemma}\label{unimin} Let $U\subset M$ be an open self-visible
set with $\Bar{U}$ compact, and let $f$ be a $C^1$ function defined on
some open neighborhood of $\Bar{U}$.

(a) If $\grad f$ is outward-pointing at each point of $\partial U$,
then $f|_{\Bar{U}}$ never achieves its minimum at a point of $\partial
U$, and hence achieves it at some critical point $\Bar{q}\in U$.

(b) If $f$ is convex on $U$ then the critical points of $f$ in $U$, if
any, are global minima of $f|_{\Bar{U}}$.  If $f$ is strictly convex
on $U$ then there is at most one critical point. \qedns

\end{lemma}

\begin{corollary}\label{onecom}
Let $Q\subset M$. Suppose that $f:M\to \bfr$ is $C^1$ on an open
neighborhood of $\Bar{\Star(Q)}$ and that for every open strongly convex
superset $U\supset Q$, the gradient of $f$ is outward-pointing along
$\partial U$. Suppose that there exists an open strongly convex
superset $U_1\supset Q$ of compact closure for which $f|_{U_1}$
achieves a minimum at some point $\Bar{q}$, and that $\Bar{q}$ is the
unique local minimum of $f$ in $U_1$. Then $\Bar{q}\in\ohull(Q)$ and
is the unique local minimum of $f$ in $\ohull(Q)$.  If ${\cal U}$ is
any collection of supersets of $Q$ on each of which $f$ has a unique
local minimum, then $\Bar{q}$ is the unique local minimum of $f$ in
$\Union_{U\in{\cal U}}U$.

\end{corollary}

\pf 
Let $U\supset Q$ be open and strongly convex, with $\Bar{U}$
compact. Then $U\intersect U_1$ is an open strongly convex superset of
$Q$, so $\na f$ is outward-pointing along $\partial (U\cap U_1)$, and
$\Bar{U\intersect U_1}$ is compact. By Lemma \ref{unimin}, $f|_{U\cap
U_1}$ achieves a minimum at some point $q$. But $\Bar{q}$ is the
unique local minimum of $f$ in $U_1$; hence $q=\Bar{q}$, so
$\Bar{q}\in U$ for every open convex superset of $Q$.
\qedns

In the case of our functions $f_Q$, the key point is that if $U$ is an
arbitrary open strongly convex superset of $Q$, then from
(\ref{defyq}) the vector field $Y_Q$ is inward-pointing along
$\partial U$, so
$\grad(f_Q)$ is outward-pointing and Corollary \ref{onecom} applies.
Thus, while strongly convex or regular geodesic {\em balls} are
essential to the proofs of Karcher's and Kendall's theorems (as well
as to the proof of Theorem
\ref{mythm} in this paper), once one has existence and uniqueness
within even one bounded strongly convex open ball, balls can
essentially be dispensed with in favor of general strongly convex open
sets.  This allows us to frame our desired characterization of {\em
the} center of mass, or {\em average}.

\begin{defn}\label{defcom2} {\rm
If $(Q,\mu)$ has a unique center of mass $\Bar{q}$ in $\ohull(Q)$, we
call $\Bar{q}$ the {\em primary} center of mass, or simply {\em the}
center of mass, of $(Q,\mu)$.  If $\tQ$ is a finite list of points and
$\mu$ is the normalized counting measure, we also refer to the primary
center of mass as \underline{the} {\em (Riemannian) average} of $\tQ$.
}
\end{defn}

Thus, combining Theorems \ref{karchthm} and \ref{kendthm} 
with Corollary \ref{onecom}, we have the following.

\begin{corollary}\label{primcom}
Suppose $Q\subset M$ is contained in a strongly convex regular
geodesic ball.  Then for any probability distribution $\mu$ on $Q$,
the primary center of mass $\Bar{q}$ of $(Q,\mu)$ exists, lies in
$\ohull(Q)$, and is the unique Karcher mean of $(Q,\mu)$ in
$\regstar(Q)$. If $f_Q$ has a local minimum at $\Bar{q}$, then the
restriction of $f_Q$ to $\regstar(Q)$ achieves its absolute minimum at
$\Bar{q}$ and nowhere else. \qedns

\end{corollary}

 In particular, Karcher means given by any two balls containing $Q$
in Karcher's or Kendall's theorem coincide.

Note that $\regstar(Q)$ can be much larger than any single regular
geodesic ball. For example, let $M$ be the unit sphere $S^n$. Let
$Q\subset S^n$ be a set of two non-antipodal points, let $C$ be the
minimal arc joining the points, and let $C_{\rm opp}$ be the arc
antipodal to $C$. Then $\regstar(Q)=\Star(Q)=S^n-C_{\rm opp}$.  In
this and some other obvious examples on spheres, $\regstar(Q)$
coincides with $IC(Q)$:= the largest open superset of $Q$ that does
not meet the cut-locus of any point of $\hull(Q)$. It is plausible
that in general $\regstar(Q)\subset IC(Q)$. However an example in
\cite{wkendall2} shows that in general $\regstar(Q)$ in Corollary
\ref{primcom} cannot be replaced by $IC(Q)$ in general without sacrificing
uniqueness.

It is plausible that Corollary \ref{primcom} remains true with
``ohull'' by ``hull'', but the author has not found a proof.  However,
Cheeger and Gromoll's general structure theorem for convex sets
(\cite{CG} Theorem 1.6; note that our ``convex'' is Cheeger and
Gromoll's ``strongly convex'') shows that $\hull(Q)$ has a
well-defined dimension. Only if this dimension equals $\dim(M)$ is our
definition of ohull exactly what is needed for the given proof of
Corollary \ref{primcom}.  However, Corollaries \ref{onecom} and
\ref{primcom} can be sharpened to include the case
$\dim(\hull(Q))<\dim(M)$; see \cite{groi1} (the original preprint
version of this paper, available from the author).


\setcounter{equation}{0}
\section{Constructing the primary center of mass}
\label{construct}

\def\deltamax{\Delta_{\rm max}}
\def\deltamin{\d_{\rm min}}

The methods of \S 2 allow us to give a constructive proof of a version
of Theorem \ref{karchthm}.  This section is devoted to the proof and a
discussion of the consequences.  Throughout we assume that the set $Q$
lies in a strongly convex ball $B$.

The vector field $Y_Q$ on $B$ gives rise to a map
$\Psi_Q=\Psi_{Y_Q}=\exp\circ Y_Q:B\to M$ as in Section \ref{zvf}.  To
apply our contracting-mapping result, Theorem \ref{quant_cor}, we need
bounds on $\|\na Y_Q + I\|$. Heuristically it is easy to understand
why this quantity is small, provided $\rho$ is small enough. Let ${\sf
g}^{-1}:T^*M\tensor T^*M\to T^*M\tensor TM\iso{\rm End}(TM)$ be the
isomorphism defined by using the metric to identify $T^*M$ with $TM$
(``raising an index'' on the second factor of $T^*M\tensor T^*M$).
For any function $f:M\to\bfr$, let $\hess(f)=\na\na f\in\Gamma({\rm
Sym}^2T^*M)$ denote its covariant Hessian, and let $\hessp(f)={\sf
g}^{-1}(\hess(f))\in\Gamma({\rm End}(TM)$.  From Theorem
\ref{karchthm}(a) we have $\na Y_Q= -\hessp(f_Q)$.  
In normal coordinates $\{x^i\}$ centered at a point $q$, for points near $q$
we have 
\ben
\hess(\frac{1}{2}r_q^2)=\sum_i dx^i\tensor dx^i\approx
\sum_{i,j}g_{ij}dx^i\tensor dx^j = g,
\een
so that $\hess'(\frac{1}{2}r_q^2)\approx {\sf
g}^{-1}g=I$ near $q$.  From \cite{karcher} Theorem 1.5 we have
\be\label{inthess}
(\na Y_Q)(p)=-\int_Q\hessp(\frac{1}{2}r_q^2)|_p\ d\mu(q).
\ee
Thus for general $Q$ contained in a small set, at points near $Q$ the
endomorphism $-\na Y_Q$ is an average of endomorphisms close to the
identity, and hence is close to the identity.

A quantitative bound on $\|\na Y_Q+I\|$ can be obtained
in terms of the functions $h_\pm,h_0$
defined by
\be\label{defhpm}
h_+(x)=x\cot x \ (0\leq x<\pi \ \mbox{only}), \ \ 
h_0(x)\ident 1, \ \ 
h_-(x)=x\coth x.
\ee
The function $h_+$ is monotone decreasing (hence $\leq 1$), while
$h_-$ is monotone increasing (hence $\geq 1$). Define
\bearray
h(\l,r)&=&h_{{\rm sign}(\l)}(|\l|^{1/2}r)
= \frac{\ch(-\l r^2)}{\sh(-\l r^2)}, \label{defh}\\
\psi(\l,r)&=& {\rm sign}(\l)(1-h(\l,r)).\label{defpsi2}
\eearray
Then $h$ is an analytic (entire) function of $\l r^2$, with
$h(\l,r)=1-\frac{1}{3}\l r^2 +O((\l r^2)^2)$. For every $\l$ the
function $r\mapsto\psi(\l,r)$ is nonnegative, monotone increasing on
$[0,\pi)$ if $\l>0$ and on $[0,\infty)$ if $\l\leq 0$, and 
$\psi(\l,r)=\frac{1}{3}|\l| r^2 +O(\l^2 r^4)$.  For 
$\d\leq\Delta \in \bfr$ and
$0\leq r<\pi\Delta^{-1/2}$ (the upper limit on $r$ applying only if
$\Delta>0$), define
\bearray\label{psimax}
\psimax(\d,\Delta,r)&=&\max(\psi(\Delta,r),\psi(\d,r))\\
&=& \frac{1}{3}|K| r^2+ O(|K|^2 r^4) \label{H2}
\eearray
where $|K|=\max(|\d|,|\Delta|)$.
Note that $\psimax$ is monotone increasing in $\Delta$ and $r$, monotone
decreasing in $\d$. Observing that $\frac{d^2}{dr^2}\psi(\pm
1,r)=\left\{\begin{array}{ll}2\csc^2r\\ 2{\rm
csch}^2r\end{array}\right\}\cdot \psi(\pm 1,r)\geq 0$, it also follows
that $\psimax$ is a convex function of each argument with the other two held
fixed.  The relevance of $\psimax$ is in the following lemma.

\begin{lemma}\label{hesslemma}
Let $p,q\in M$ with $d(p,q)<\rinj(q)$ and let $\d$ and $\Delta$ be
lower and upper bounds, respectively, for the sectional curvatures of
$M$ along the minimal geodesic from $q$ to $p$ $\g$; if $\Delta>0$
also assume $d(p,q)<\pi\Delta^{-1/2}$. Then
\be\label{hessbound2} 
\|\hessp(\frac{1}{2}r_q^2)-I\|(p)
\leq \psimax(\d,\Delta,d(q,p)).
\ee
If $d(p,q)\cdot\max(0,\Delta)^{1/2}<\pi/2$, then 
\be\label{hessbounds}
\hess(\frac{1}{2}r_q^2)|_p>0.
\ee
\end{lemma}

\pf Both statements follow immediately from Lemma \ref{hesslemma_app} in 
the Appendix. \qedns

Henceforth we assume that $Q$ lies in a ball $B_D(p_0)$ and analyze
the vector field $Y_Q$ on a possibly larger concentric ball
$B=B_\rho(p_0)$, still assumed strongly convex.  We apply the lemma to
points $p\in B,q\in Q$, setting $\d=\d(B), \Delta=\Delta(B)$. For such
points we have $d(p,q)<\rho+D$, so to meet the potential restriction
on $d(p,q)$ in the lemma, we assume that
$(\rho+D)\max(0,\Delta)^{1/2}<\pi$.  From (\ref{inthess}),
(\ref{hessbound2}), and the monotonicity of $\psimax$ we then have
\be\label{bddbyh}
\|\na Y_Q + I\|=\|\int_Q(\hessp(\frac{1}{2}r_q^2)-I)\ d\mu(q)\| \leq
\psimax(\d,\Delta,\rho+D). 
\ee
We also have 
\be\label{supbound}
\|Y_Q(p)\|\leq \sup_{q\in Q}d(p,q)\leq \rho+D.
\ee
Hence from Theorem \ref{quant_cor}, for all $p\in B$ we have
\be\label{kq}
\|(\Psi_Q)_{*p}\|\leq 
\k(p_0;\rho,D) := 
\phi_\pm((\rho+D)\aku^{1/2})+C_1(\d,\rho+D)\ \psimax(\d,\Delta,\rho+D).
\ee
where $\aku=\aku(B_{\rho}(p_0))$, and where the choice of sign in
$\phi_\pm$ is governed by the following convention.

\bsn{\bf Notation Convention 4.2}\addtocounter{theorem}{1}
{\em For the remainder of this paper, when an expression of the form
$\phi_\pm(x)$ appears, $\phi_+(x)$ is to be used if $M$ is a locally
symmetric space of nonnegative curvature and $x\leq 3\pi/4$;
$\phi_-(x)$ is to be used otherwise. }

\bs 
To ensure that $\Psi_Q$ is a contraction we want $\k(p_0;\rho,D)<1$,
which will be true for small $\rho$ since $\phi_\pm(x)$ and
$\psimax(\cdot,\cdot,x)$ are $O(x^2)$. This is not enough by itself to
ensure existence of a fixed point:

\begin{defn}\label{deftether}{\rm  
Call a map $\Psi: ({\rm domain}(\Psi)\subset M)\to M$ {\em tethered to
$Q$} if, for every strongly convex regular geodesic ball $B$
containing $Q$, (i) $\Psi$ is defined on $B$ and (ii) $\Psi(B)\subset
B$.
 }
\end{defn}

If we knew $\Psi_Q$ to be tethered to $Q$ (which implicitly requires 
${\rm domain}(\Psi)\supset \regstar Q$), we could apply the general
form of the Contracting Mapping Theorem (which assumes {\em a priori}
that the contracting map preserves its domain) to conclude that
$\Psi_Q$ has a unique fixed point in $\Bar{B_\rho(p_0)}$ as long as
$\k(p_0;\rho,D)<1.$ In Euclidean space, $\Psi_Q$ is always tethered to
$Q$ trivially: $\Psi_Q$ maps the entire space to a single point
contained in the convex hull of $Q$.  On a general manifold, if $Q$
consists of a single point then $\Psi_Q$ is tethered to $Q$ for the same
trivial reason. Thus it seems likely that on general $M$, tethering
will occur provided $\diam(Q)$ is sufficiently small.  It is plausible
that this happens for any $Q$ contained in a strongly convex regular
geodesic ball, but the author has neither a proof nor a
counterexample.  The lack of such a proof is the sole reason that in
our center-of-mass application we use Theorem
\ref{cmt} (in the guise of Theorem \ref{quant_cor}) rather than the
more general Contracting Mapping Theorem (but note that Theorem
\ref{quant_cor} may still be needed in other applications, i.e. those using 
maps $\Psi_Y$ with $Y$ not of the form $Y_Q$, since most such general
maps will not be tethered).  The cost is that the upper bound on the
diameter of $Q$ (or other measures of size such as the
``circumradius'') for which we can ensure that $\Psi_Q$ has a fixed
point is smaller than it would be if we knew that tethering
occurred. Since it may be possible to prove tethering, either in
general or in specific cases, in the remaining theorems of this paper
we include statements of what one can conclude in the tethered case.

Assuming $\k(p_0;\rho,D)<1$, to conclude from Theorem
\ref{quant_cor} that $\Psi_Q$ has a fixed point, we additionally need
to have
\be\label{needy}
\|Y_Q(p_0)\|< (1-\k(p_0;\rho,D))\rho :=s(p_0;\rho,D).
\ee
Clearly (\ref{supbound}) is of no help here.  However, the left-hand
side of (\ref{needy}) does not depend intrinsically upon $\rho$, but
only upon $(Q,\mu)$. We are taking $\rho\geq D$, so furthermore
$s(p_0;\rho,D)\geq s(p_0;\rho,\rho):=s_2(p_0;\rho)$. The basis
of the argument over the next few pages is simply that as long as
$\|Y_Q(p_0)\|$ is less than the maximum value of the function
$s_2(p_0; \cdot)$, there will be {\em some} radius $\rho$ for which
(\ref{needy}) is satisfied even with $D=\rho$, hence for all
$D\leq\rho$ as well.  

Note also that $\|Y_Q(p_0)\|\leq D$, so that an upper bound on $D$
implies an upper bound on $\|Y_Q(p_0)\|$.  Thus the most general
conclusions we eventually draw will be those that have an upper bound
only on $\|Y_Q(p_0)\|$ (hence on $(Q,\mu)$) as a hypothesis, but as a
corollary all such conclusions hold with an upper bound on $D$, a more
easily checked and therefore more practical hypothesis. Eventually in
Corollary \ref{easycor} we will take $p_0$ to lie in $Q$, which will
give us even more control since we can then take $D={\rm diam}(Q)$.

Since we are interested not just in the {\em existence} of ``good''
radii $\rho$ and $D$, but on estimating their size, we first prove a
lemma establishing some properties of the function $s$; these will be
used to estimate the size of balls on which $\Psi_Q$ has a fixed
point.  In practice one is usually not presented with an explicit
growth rate for $|\d|,|\Delta|,$ or $|K|$ as functions of $\rho$ in
(\ref{kq}), so we also examine the consequences of a (potentially
less sharp but usually more practical version of the bound in
(\ref{needy}), replacing the function $s$ by a function $\tilde{s}$
defined below. The sharp bounds, however, are needed for the best
estimates in \cite{groi2} for an averaging algorithm on size-and-shape
spaces.

\begin{defn}\label{def44}{\rm 
 Let $p\in M$. 
(a) For $0\leq D\leq\rho<\rreg(p)$, let
$\Delta_{p,\rho}=\Delta(B_\rho(p))$, $\d_{p,\rho}=\d(B_\rho(p))$,
$\aku_{p,\rho}=\aku(B_\rho(p))$, and

\bearray
\label{defkr}
\k(p;\rho,D)&=&\phi_\pm ((\rho+D) |K|_{p,\rho}^{1/2})+C_1(\d_{p,\rho},\rho+D)
\psimax(\d_{p,\rho},\Delta_{p,\rho},\rho+D),\\
\label{defsr}
s(p;\rho,D)&=& (1-\k(p;\rho,D))\rho.
\eearray
(If $\d_\rho=-\infty$ interpret (\ref{defkr}) as
$\k(p;\rho,D)=\infty$.)  

(b) Let $r_1\in (0,\rreg(p))$, and let $\tilde{\Delta}(\cdot)$
(respectively $\tilde{\d}(\cdot)$) be any continuous monotonically
increasing (resp. decreasing) function on $[0,r_1]$ such that
$\Delta_{p,\rho}\leq\tilde{\Delta}(\rho)$, $\d_{p,\rho}\geq
\tilde{\d}(\rho)$,
$r_1\cdot\max(0,\tilde{\Delta}(r_1))^{1/2}<\pi/2$. For $0\leq
D\leq\rho\leq r_1$ define $\tilde{\k}(p,\tilde{\Delta},\tilde{\d};\rho,D)$
to be the right-hand side of (\ref{defkr}) with $\Delta_{p,\rho},
\d_{p,\rho}, |K|_{p,\rho}$ replaced by $\tilde{\Delta}(\rho),
\tilde{\d}(\rho)$, $\max(|\tilde{\Delta}(\rho)|, |\tilde{\d}(\rho)|)$
respectively, and define
\be
\label{defts}
\tilde{s}(\rho,D)=\tilde{s}(\tilde{\Delta},\tilde{\d};\rho,D)=
(1-\tilde{\k}(\tilde{\Delta},\tilde{\d};\rho,D))\rho. 
\ee
}
\end{defn}

In practice, $\tilde{\Delta}$ and $\tilde{\d}$
will usually be {\em constant} functions, global upper and lower
curvature bounds on $B_{r_1}(p)$. We define $\tilde{s}$ in greater
generality above because this enables not only stronger results, but
shorter proofs: anything proven for the more general functions
$\tilde{s}$ applies to the special case $\tilde{s}=s$.

We construct from such a function $\tilde{s}$ several numbers and
functions of $D$: $\tilde{D}_{\rm crit}, \tilde{D}_{\rm max}$, and
$\tilde{\rho}_i$, all defined below. The meaning of the
$\tilde{\rho}_i(p,r_1;D)$ is indicated by the $\rho_i$ in Figure 1; the
qualitative correctness of Figure 1 is proven in Lemma \ref{slemma}.

\vbox{\medskip
\centerline{
\includegraphics[width=3in]{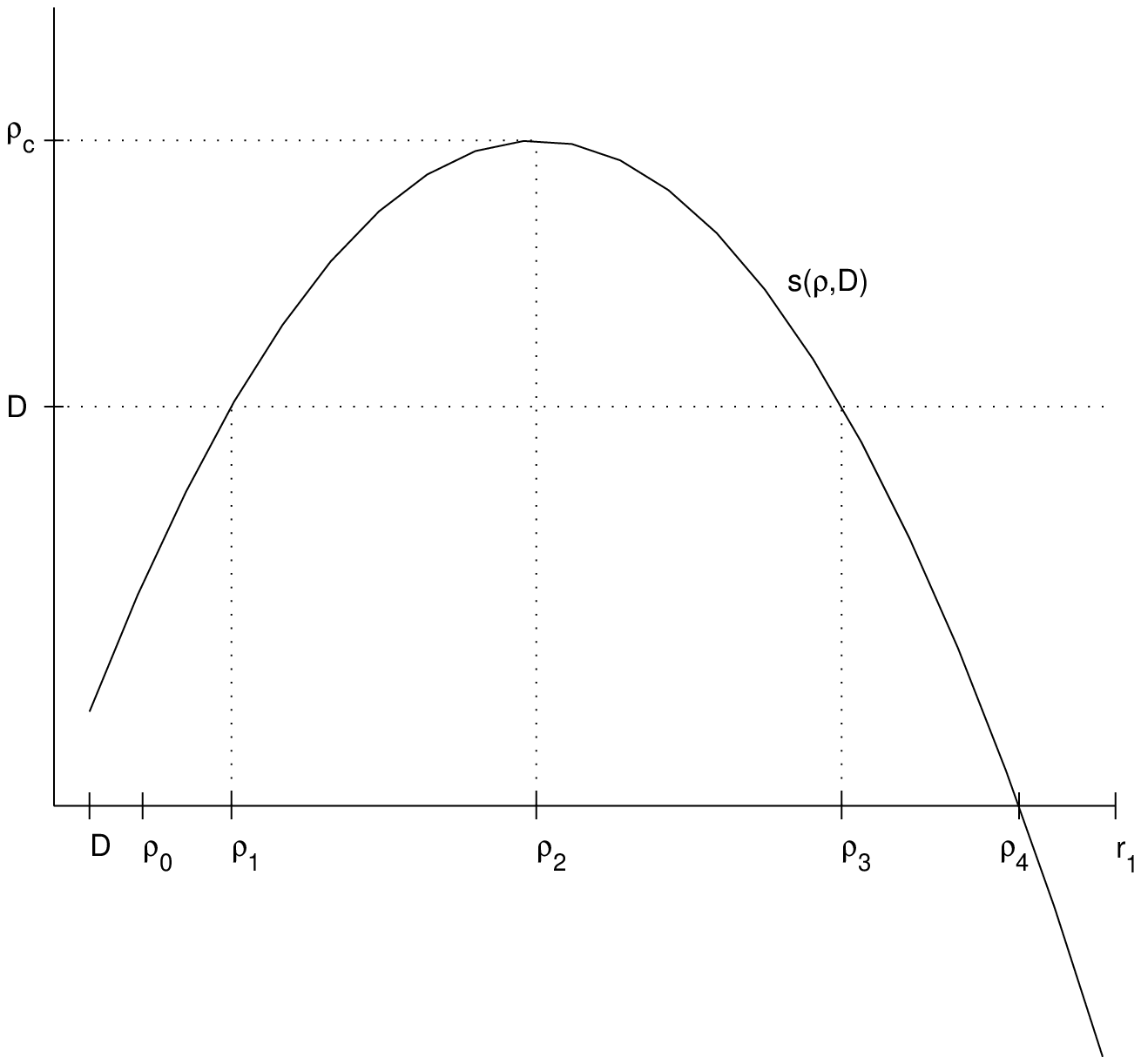}
}\medskip
\footnotesize 
\ni {\bf Figure 1.} A (not-to-scale) 
sketch of $s(\rho,D)$ versus $\rho$ for some fixed $D<D_{\rm crit}$,
assuming $r_1>\rho_4$.  A sketch of $\tilde{s}(\rho,D)$ for a fixed
$D<\tilde{D}_{\rm crit}$ would be similar, with the $\rho_i$ replaced
by $\tilde{\rho}_i, 1\leq i\leq 4$.  $D_{\rm crit}$ is the maximum
value of $s(\rho,D_{\rm crit})$; for $D>D_{\rm crit}$, the graph of
$s(\rho,D)$ lies entirely below the horizontal line at height $D$.  To
illustrate the maximum number of distinct radii we have sketched the
case in which $\rho_4$ is strictly less than $r_1$, i.e. in which
$\k(\rho,D)$ reaches 1 before $\rho$ reaches $r_1$.  The picture for
smaller $r_1$ can be obtained from this one by moving $r_1$ to the
left, say to $r_{1,{\rm new}}$, truncating the diagram to the right of
$r_{1,{\rm new}}$ and decreasing $D$, if necessary, to keep it less
than the maximum value of $s$ on $[D,r_{1,{\rm new}}]$ (hence keeping
$r_{1,{\rm new}}>\rho_1(D_{\rm new})$). If any of $\rho_2, \rho_3,
\rho_4$ in the picture above is to the right of $r_{1,{\rm new}}$, the
corresponding $\rho_{i,{\rm new}}$ is defined to be $r_{1,{\rm new}}$.
}
\bigskip

The definitions of $\tilde{D}_{\rm crit}(p,r_1), \tilde{D}_{\rm
max}(p,r_1)$ and the $\tilde{\rho}_i(p,r_1; \cdot)$ are given in
(\ref{defdmax}--\ref{defdcrit}) and (\ref{defrho0t}--\ref{defrho1t})
below.  Here and below we suppress the parameters $\tilde{\Delta}$ and
$\tilde{\d}$ rather than write $\tilde{D}_{\rm
crit}(p,r_1,\tilde{\Delta}, \tilde{\d})$ etc.; these parameters are
always present implicitly.  For the sharp-curvature-bound case
($\tilde{s}=s$), we omit the tildes and just write ${D}_{\rm
crit}(p,r_1), {D}_{\rm max}(p,r_1)$ and ${\rho}_i(p,r_1)$. Since
$\tilde{\k}(\cdot,\cdot,0,0)=0$, the sets over which the suprema are
taken below are nonempty and the suprema well-defined.

\bearray
\tilde{D}_{\rm max}(p,r_1)
&=&\sup\{D\in [0,r_1] \mid
\tilde{\k}(D,D)<1\}. \label{defdmax}
\\
D_{\rm max}(p)
&=&\sup\{D\in [0,\rreg(p)) \mid
\tilde{\k}(D,D)<1\} \label{defdmax2}\\
\nonumber&=& \sup\{{D}_{\rm max}(p,r_1)\mid r_1<\rreg(p)\}.
\\ \nonumber \\
\label{defdcritt}\label{defdcrit}
\tilde{D}_{\rm crit}(p,r_1)&=&
\sup\{D\in[0,r_1]\mid \exists \rho\in[D,r_1]\ \mbox{for which} \
\tilde{s}(\rho,D)>D \}. \\ \nonumber
D_{\rm crit}(p)
&=&\sup\{D\in [0,\rreg(p)) 
\mid \exists \rho\in[D,\rreg(p))\ \mbox{for which} \
\tilde{s}(\rho,D)>D \} \\ \label{defdcrit2} \\
&=& \sup\{D_{\rm crit}(p,r_1)\mid r_1<\rreg(p)\}.\nonumber
\eearray

\be
\label{defrho0t}
\tilde{\rho}_0(p,r_1;D)=\tilde{\rho}_0(p;D)
=\left\{ \begin{array}{ll} D/h_+(2D\tilde{\Delta}(D)^{1/2}) &
{\rm if}\ \tilde{\Delta}(D)>0,\\
D & {\rm if}\ \tilde{\Delta}(D)\leq 0.
\end{array}\right.
\ee

\ni 
For $0\leq D <\tilde{D}_{\rm max}(p,r_1)$, define
\bearray
\tilde{\rho}_4(p,r_1;D)&=&
\sup\{\rho\in[D,r_1]\mid \tilde{\k}(\rho,D)<1\};
  \label{defrho4t}
\eearray
\ni for  $0\leq D <{D}_{\rm max}(p)$ define
\bearray
{\rho}_4(p;D)&=&
\sup\{\rho\in[D,\rreg(p)\mid \tilde{\k}(\rho,D)<1\}.  \label{defrho4}
\eearray

\ni For $0\leq D<\tilde{D}_{\rm crit}$ define
\bearray
\label{defrho3t}
\tilde\rho_3(p,r_1;D)&=&\sup\{\rho\in [0,r_1] \mid
\tilde{s}(\rho,D)>D \}, \\
\label{defrho1t}
\tilde\rho_1(p,r_1;D)&=&\inf\{\rho\in [0,r_1] \mid
\tilde{s}(\rho,D)>D \};
\eearray
\ni
for  $0\leq D<{D}_{\rm crit}$ define
\bearray
\label{defrho3}
\rho_3(p;D)&=&\sup\{\rho\in [0,\rreg(p)) \mid
\tilde{s}(\rho,D)>D \}, \\
\label{defrho1}
\rho_1(p;D)&=&\inf\{\rho\in [0,\rreg(p))  \mid
\tilde{s}(\rho,D)>D \};
\eearray

Note that for $i=1,3,4,$ $\rho_i(p;D)$ can alternatively be written as
a supremum or infimum (over $r_1$) of $\rho_i(p,r_1;D)$ as we did
above for $D_{\rm max}(p)$ and $D_{\rm crit}(p)$. Note also that in
(\ref{defdmax}) and (\ref{defrho4t}), ``$\tilde{\k}(\cdot,\cdot)<1$''
can be replaced by ``$\tilde{s}(\cdot,\cdot)>0$'' without altering the
definitions of $\tilde{D}_{\rm max}$ and $\tilde{\rho}_4$.

The technical lemma below establishes some useful properties of the
objects just defined, including monotonicity in parameters.

\ms
\begin{lemma}\label{slemma}
Let $p\in M$ and let $r_1\in(0,\rreg(p))$.  Let
$\tilde{\Delta},\tilde{\d}$ be continuous monotone bounds on curvature
as in Definition \ref{def44}(b), and let $\tilde{D}_{\rm crit}=
\tilde{D}_{\rm crit}(p,r_1),
\tilde{D}_{\rm max}=\tilde{D}_{\rm max}(p,r_1)$, and 
$\tilde{\rho}_i(\cdot)=\tilde{\rho}_i(p,r_1;\cdot)$ be as in
(\ref{defdmax}--\ref{defrho1t}).

For $D\in [0,r_1]$ let ${J}_D=\{\rho\in[0,r_1] \mid
D<\tilde{s}(\rho,D)\}$. For each $D$, the set ${J}_D$ is either empty
or an interval with endpoints $\tilde{\rho}_1(D),\tilde{\rho}_3(D)$.
If $D_2> D_1$ then $\Bar{{J}_{D_2}}\subset {J}_{D_1}$, so $\{D\mid
J_D\neq\emptyset\}$ is an interval whose right endpoint is
$\tilde{D}_{\rm crit}$. $\tilde{D}_{\rm crit}>0$ and $\bigcap_{0\leq
D< \tilde{D}_{\rm crit}}J_D$ consists of a single point
$\tilde{\rho}_{\rm crit}$, satisfying $\tilde{D}_{\rm
crit}=\tilde{s}(\tilde{\rho}_{\rm crit},\tilde{D}_{\rm
crit})=\max_{\rho\in[D,r_1]}s(\rho,\tilde{D}_{\rm crit}),$ the maximum
being achieved uniquely. The following are true:
\begin{enumerate}
\item
$\tilde{D}_{\rm max}\geq\tilde{D}_{\rm crit}$, with equality if and only
if $\tilde{D}_{\rm crit}=r_1$.

\item
$D_{\rm crit}\geq \tilde{D}_{\rm crit}, D_{\rm max}\geq\tilde{D}_{\rm
max}$.

\item
$\rho_4(D)\geq\tilde{\rho}_4(D)$ for all $D<\tilde{D}_{\rm max}$.

\item \label{dmaxbound} 
$(2\tilde{D}_{\rm max})\cdot\max(0,\tilde{\Delta}
(\tilde{D}_{\rm max}))^{1/2}<\pi/2.$

\item \label{kappaconvex}
For each $D\in [0,r_1]$, the function $\rho\mapsto \tilde{\k}(\rho,D)$ 
on $[D,r_1]$ is continuous, monotone
increasing, and convex. The function $\rho\mapsto
\tilde{\k}(\rho,\rho)$ is $O(|\tilde{K}|_\rho\ \rho^2),$ where
$|\tilde{K}|_\rho =\max(|\tilde{\Delta}(\rho)|, |\tilde{\d}(\rho)|)$.

\item \label{sconcave}
For each $D\in [0,r_1]$,
the function $\tilde{s}(\cdot,D)$ is concave and achieves its maximum
at a unique point $\tilde{\rho}_2(D)\in (0,r_1]$.

\end{enumerate}

\ni 
For each $D<\tilde{D}_{\rm crit}$ the following are true, where
$\tilde{\rho}_i=\tilde{\rho}_i(D).$ 

\begin{enumerate}\setcounter{enumi}{6}
\item\label{conc7}
$\rho_3\geq\tilde{\rho}_3,\ \rho_0\leq\tilde{\rho}_0,$ and
$\rho_1\leq\tilde{\rho}_1$.

\item \label{it0}\label{it1}
The following order-relations hold (cf. Figure 1):
\bearray
&D\leq
\tilde{\rho}_0\leq\tilde{\rho}_1<\tilde{\rho}_{\rm crit} 
\leq\tilde{\rho}_3<\tilde{\rho}_4\leq \tilde{D}_{\rm max}\leq r_1.&
\label{rho24}\label{rho01}\label{rho1crit3}
\eearray

\item \label{rho4bound} 
$(\tilde{\rho}_4+D)\cdot\max(0,\tilde{\Delta}
(\tilde{\rho}_4))^{1/2}<\pi/2.$

\end{enumerate}
As a special case, all conclusions above are true with the tildes
erased. As a corollary, conclusion \ref{dmaxbound} is true also with
$\tilde{D}_{\rm max}(p;r_1)$ replaced by $\tilde{D}_{\rm max}(p)$;
conclusions \ref{kappaconvex} and \ref{sconcave} are true with the
tildes erased and with $[D,r_1]$ replaced by $[D,\rreg(p))$; and
conclusions \ref{conc7}--\ref{rho4bound} are true with $\tilde{D}_{\rm
crit}(p;r_1)$ and $\tilde{D}_{\rm max}(p,r_1)$ replaced by 
${D}_{\rm crit}(p)$ and ${D}_{\rm max}(p)$ respectively, 
$\tilde{\rho}_i(p,r_1;D)$ replaced by $\rho_i(p;D)$, and
``$\tilde{\rho}_4\leq r_1$'' replaced by ``$\rho_4<\rreg(p)$''.

\end{lemma}

\pf
From the definition of $\tilde{\k}$ continuity in all parameters is
clear, and it is easy to check that $\tilde{\k}(\rho,D)\leq
\tilde{\k}(\rho,\rho)=O(|\tilde{K}|_\rho\
\rho^2)$.  We have already noted that $\psimax(\d,\Delta,r)$ is
monotone increasing in $r$ and $\Delta$, decreasing in $\d$, and
convex in each variable separately; the same is true of
$C_1(\d,r)$. The functions $\phi_\pm$ are monotone increasing and
convex. Monotonicity and convexity of $\phi_\pm,H,$ and $C_1$ are
retained after composition with the monotone functions
$\tilde{\d}(\rho), \tilde{\Delta}$.

It follows that with $D$ held fixed, $\tilde{\k}(\cdot, D)$ is
continuous, monotone increasing and convex, and hence that
$\tilde{s}(\cdot, D)$ is continuous, concave, and, because of the
factor of $\rho$ in (\ref{defts}) and monotonicity, nonconstant on any
interval of positive length. Since $\tilde{\k}(\rho,0)=O(\rho^2)$,
$\tilde{s}(\rho,0)>0$ for $\rho>0$ sufficiently small. Hence
$J_0$ is nonempty, and by continuity so is $J_D$ for
sufficiently small positive $D$. Hence $\tilde{D}_{\rm crit}>0$.

For each fixed $D$, the concavity and local nonconstancy of the
function $\tilde{s}(\cdot,D)$ implies that its maximum value
$\tilde{\rho}_c(D)$ on $[0,r_1]$ is achieved at a unique point
$\tilde{\rho}_2(D)$, and for any $a<\tilde{\rho}_c(D)$
the set $\{\rho\in[0,r_1]\mid
\tilde{s}(\rho)> a\}$ is an interval; 
in particular each set $J_D$ is an
interval. Since $D_2>D_1$ implies $\tilde{s}(\rho,D_2)<
\tilde{s}(\rho,D_1)$ strictly for $\rho>0$, the asserted nesting of 
the intervals $J_D$ also follows. The intersection of the nonempty
$J_D$ is nonempty because their closures are nested, and the
intersection has only one point $\tilde{\rho}_{\rm crit}$ since
$\tilde{s}(\cdot,D_{\rm crit})$ is nowhere constant.  Continuity
implies $\tilde{D}_{\rm crit}=\tilde{s}(\tilde{\rho}_{\rm crit},
\tilde{D}_{\rm crit})$.

From its definition clearly $\tilde{s}(\rho,D)\leq \rho$.  All the
inequalities asserted in statement \ref{it0} follow immediately from
the foregoing, except for $\tilde{\rho}_0\leq \tilde{\rho}_1$. The
latter inequality follows from chasing through the definitions and
monotonicity of the ingredients in $\tilde{\k}$. A helpful observation
is that from (\ref{defkr}) we have
\be\label{hint}
\tilde{\k}(\rho,D) \geq \psimax(\tilde{\d}(\rho),
\tilde{\Delta}(\rho), \rho+D)
\geq \psi(\tilde{\Delta}(\rho), \rho+D).
\ee
It also follows that $\tilde{s}(\tilde{D}_{\rm crit}, \tilde{D}_{\rm
crit})
\geq \tilde{s}(\tilde{\rho}_{\rm crit}, \tilde{D}_{\rm crit})= 
\tilde{D}_{\rm crit}>0$, so that $\tilde{D}_{\rm max}\geq\tilde{D}_{\rm
crit}$. 

The monotonicity of $\phi_\pm,C_1,$ and $\psimax$ imply that if
$\rho\leq r_1$, then $\k(\rho,D)\leq
\tilde{\k}(\rho,D)$, and hence $\tilde{s}(\rho,D)\leq s(\rho,D)$.
Hence $\tilde{D}_{\rm crit}\leq D_{\rm crit}, \
\tilde{\rho}_1\geq\rho_1$, 
and $\tilde{\rho}_i\leq\rho_i$ for $i=3,4$.

To establish statements \ref{dmaxbound} and \ref{rho4bound} we claim first
that for $D<D_{\rm crit}$ we have
\be\label{hint2}
(\tilde{\rho}_1(D)+D)\max(0,\tilde{\Delta}(\tilde{\rho}_1(D)))^{1/2}<\pi/2.
\ee
This is true for $D=0$, so if it is false for some $D<D_{\rm crit}$
then there exists $D\in(0,D_{\rm crit})$ for which
$\tilde{\Delta}(\tilde{\rho}_1(D))>0$ and
$(\tilde{\rho}_1(D)+D)\tilde{\Delta}(\tilde{\rho}_1(D))^{1/2}=\pi/2$,
the latter implying $\psi(\tilde{\Delta}(\tilde{\rho}_1(D)),
\tilde{\rho}_1(D)+D)=1$.  But the combination $D>0, \tilde{\Delta}>0$
implies strict inequality in (\ref{hint}), so
$\tilde{\k}(\tilde{\rho}_1(D),D)>1$ and
$\tilde{s}(\tilde{\rho}_1(D),D)<0$; but from the definition of
$\tilde{\rho}_1$ we have $\tilde{s}(\tilde{\rho}_1(D),D)\geq D$.
Hence (\ref{hint2}) holds for all $D<D_{\rm crit}$.  Therefore if
statement \ref{rho4bound} is false, there exists
$\rho\in(\tilde{\rho}_1(D), \tilde{\rho}_4(D)$ for which
$(\rho+D)\tilde{\Delta}(\rho)^{1/2}=\pi/2.$ From (\ref{hint}) we again
conclude that $\tilde{\k}(\rho,D)>1$, and since $\rho\geq \rho_1(D)>0$
this implies the strict inequality $\tilde{s}(\rho,D)< 0,$ a
contradiction since $\rho\in (0,\tilde{\rho}_4(D))$.  This proves
statement \ref{rho4bound}; a shorter version of the same argument
yields statement \ref{dmaxbound}. \qedns

\begin{remark}{\rm
In Definition \ref{def44} and Lemma \ref{slemma}, the restriction
``$r_1< \rreg(p)$'' can be replaced by the less restrictive
``$r_1\cdot\max(0,\Delta_{p,r_1})^{1/2}<\pi/2$''.
}
\end{remark}

\begin{corollary} \label{weakkarch}
Let $p_0\in M$, $0<r_1<\rregcvx(p_0)$. Let $D_{\rm crit},D_{\rm
max}, \rho_4,\rho_1$ be as in (\ref{defdmax2})--(\ref{defrho1}). For
$0<\rho\leq r_1$ write $B_\rho$ for $B_\rho(p_0)$. Let $Q\subset
B_{\rho_4}$ be equipped with a probability measure $\mu$, and define
$Y_Q$ and $f_Q$ by (\ref{defyq}--\ref{deffq}).  Then $Y_Q$ has at most
one zero in $B_{\rho_4}$ (equivalently, $f_Q$ has at most one critical
point in this ball); at such a zero $f_Q$ achieves its minimum value
on $B_{\rho_4}$ (in fact, on $\regstar(Q)$). If $D<D_{\rm crit}$ and
$Q\subset \Bar{B_D}$ (or more generally if $\|Y_Q(p_0)\|\leq D$), then
$Y_Q$ has a unique zero $\Bar{q}$ in $B_{\rho_4}$, and $\Bar{q}$ lies
in $\Bar{B_{\rho_1}}.$ Hence $(Q,\mu)$ has at most one center of mass
in $B_{\rho_4}$, and has exactly one center of mass in $B_{\rho_4}$ if
$\Bar{Q}\subset B_{D_{\rm crit}}$.  If $\Psi_Q$ is tethered to $Q$,
these conclusions hold with $D_{\rm crit}$ replaced by the (never
smaller and usually larger) number $D_{\rm max}$.
\end{corollary}

We will prove this simultaneously with Theorem \ref{mythm} below.  But
first, taking $r_1$ close to $\rregcvx(p_0)$ in Corollary
\ref{weakkarch}, note that Lemma
\ref{slemma} implies that the restriction on 
the radius of the ball containing $Q$ in Corollary
\ref{weakkarch} is more stringent than in Theorem \ref{karchthm}(c).  
Similarly Lemma \ref{slemma} implies that the conclusion
$\Bar{q}\in\Bar{B_{\rho_1}}$ in the corollary above is not as sharp as
Karcher's conclusion $\Bar{q}\in
\Bar{B_{\rho_0}}$, and that
the conclusion above concerning existence of at most one center of
mass in $B_{\rho_4}$ is weaker than Kendall's conclusion---at most one
center of mass in $B_{\rreg(p_0)}$---which is itself weaker than the
uniqueness and minimization statement in Corollary \ref{primcom}.
(However, we will see in \S
\ref{examples} that if $(M,g)$ has non-negative curvature,
then for $D<\tilde{D}_{\rm crit}$ the uniqueness statement in
Corollary \ref{weakkarch} is actually stronger than Karcher's.) In
fact, in view of Corollary \ref{primcom}, $B_{\rho_4}$ can be replaced
by $\regstar(Q)$ in the conclusions (but not the hypotheses) of
Corollary \ref{weakkarch}.

Thus, were Corollary \ref{weakkarch} the only outcome of the
contracting-mapping approach, we would have gained little from
it. However, the contracting-mapping approach additionally provides an
{\em algorithmic construction} of the center of mass, one that is
easily implemented in spaces for which the exponential map and its
inverse are explicitly known, and in particular for shape spaces. In
practice, any algorithm intended to average a list $Q$ of points in a
space is initialized at a point $q_0\in Q$, but there are questions of
whether the algorithm converges and whether its limit (if any) depends
on the choice of initial point.  As mentioned in the introduction, GPA
algorithms converge quite rapidly in practical applications, but it is
not readily apparent why this happens.  For a given algorithm, one may
be able to prove initial-point independence of the limit by one
argument, and convergence by another, and perhaps estimate the
convergence rate still another way. However, the contracting-mapping
approach allows one to answer all these questions at once (although
answering them individually by other means may lead to sharper
answers, as in \cite{Le1} Proposition 3, for initial-point independence in
the GPA-S algorithm).  Thus the added value of this approach lies in
the following theorem, in which we state only those direct conclusions
of the contracting-mapping approach neither contained in nor relying
on Karcher's and Kendall's theorems (except for the use of $\rho_0$ in
conclusion \ref{newit3}). In
\S \ref{conv_rate} we will see that by estimating the convergence rate of
$\Psi_Q^n(p_0)$ and combining this with Kendall's uniqueness result,
we can considerably strengthen certain parts of Theorem \ref{mythm};
see Theorem \ref{strongerthm2}.  In statement \ref{newit3} of the theorem
below, note that with the indicated restrictions on $D$, existence of
the primary center of mass is guaranteed by Corollary \ref{weakkarch},
as well as by Theorem \ref{karchthm}.

\begin{theorem}\label{mythm}
Let $p_0\in M$, $0<r_1<\rregcvx(p_0)$; for $0<\rho\leq r_1$ write
$B_\rho$ for $B_\rho(p_0)$. Let $\tilde{\Delta}(\cdot),
\tilde{\d}(\cdot)$ be continuous monotone upper and lower bounds on
curvature as in Definition \ref{def44}(b).  Let $Q\subset B_{r_1}$ be
equipped with a probability measure $\mu$, and define $Y_Q,f_Q$ by
(\ref{defyq}--\ref{deffq}).  Then, using the notation
(\ref{defdmax})--(\ref{defrho1}) with the parameter $p_0$ suppressed,
the following are true.

\begin{enumerate}
\item\label{newit1}
$\tilde{D}_{\rm max}(r_1)\leq D_{\rm max}(r_1)\leq D_{\rm max}$
and $\tilde{D}_{\rm crit}(r_1)\leq D_{\rm crit}(r_1)\leq D_{\rm
crit}$. In particular if $D<\tilde{D}_{\rm crit}(r_1)$ then all the
$\rho_i(D)$ are defined, and 
\be
D\leq \tilde{\rho}_0(D)\leq\tilde{\rho}_1(D)<\tilde{\rho}_{\rm crit}
< \tilde{\rho}_3(D)\leq\tilde{\rho}_4(D)\leq\tilde{D}_{\rm max}(r_1)
\leq r_1 \label{drho3}
\ee
where $\tilde{\rho}_{\rm crit}$ is value of $\rho$ that maximizes
$\tilde{s}(\rho,\tilde{D}_{\rm crit})$.

\item\label{it2}
For all $D\in (0,r_1]$, if  $Q\subset B_D$ and 
$\rho<\tilde{\rho}_4(D)$, then the map
$\Psi_Q=\exp\circ Y_Q:B_{\rho}\to M$ is a contraction with constant
$\tilde{\k}(p_0;\rho,D)$.

\item\label{newit3}
Assume that $Q\subset \Bar{B_D}$ (or more generally that
$\|Y_Q(p_0)\|\leq D$) and that either 

\ss\hspace{.3in}
(i) $D<\tilde{D}_{\rm crit}$ and ${\rho}_1(D)<\rho<{\rho}_3(D)$, or

 \hspace{.3in} (ii) $D<\tilde{D}_{\rm max}$, $\Psi_Q$ is tethered to
$Q$ (Definition \ref{deftether}), and $D\leq \rho<\rho_4(D)$.

\ss Then $\Psi_Q$ preserves preserves each ball $B_\rho$. In
particular this holds for the $D$-independent radius
$\tilde{\rho}_{\rm crit}$.  The sequence
of iterates $\Psi_Q^n(q)$ converges to the primary center of mass
$\Bar{q}$ of $(Q,\mu)$ for every $q\in B_{{\rho}_3(D)}$ if (i) holds,
and for every $q\in B_{{\rho}_4(D)}$ if (ii) holds. In either case
$\Bar{q}$ lies in $\overline{B_{\rho_0(D)}}\cap\ohull(Q)$.

\item\label{newit4} For $D<\tilde{D}_{\rm crit}$ the following 
relations hold:
\be
\rho_0(D) \leq \tilde{\rho}_0(D),  \ \ \ 
\rho_1(D) \leq \tilde{\rho}_1(D),  \ \ \ 
\rho_3(D)\geq \tilde{\rho}_3(D), \ \ \ \rho_4(D)\geq \tilde{\rho}_4(D).
\ee
If the curvature bounds $\tilde{\Delta},\tilde{\d}$ are taken to be
constants (e.g. $\tilde{\Delta}\ident\Delta(B_{r_1}), \tilde{\d}\ident
\d(B_{r_1})$), then 
the lower bound $\tilde{D}_{\rm crit}$ on $D_{\rm crit}$ is a
universal function of the numbers $r_1,\tilde{\Delta}$, and
$\tilde{\d}$, depending in no other way on the geometry of
$(M,g)$. Similarly the lower bounds $\tilde{D}_{\rm max}$ on $D_{\rm
max}$, $\tilde{\rho}_i$ on $\rho_i$ for $3\leq i\leq 4$, and the upper
bounds $\tilde{\rho}_i$ on $\rho_i$ for $0\leq i\leq 1$, are universal
functions of $r_1,\tilde{\Delta}$, $\tilde{\d}$, and $D$.

\end{enumerate}

\end{theorem}

\begin{remark}{\rm
The chief point of the last two sentences in Statement
\ref{newit4} is that $D_{\rm crit}$, the critical upper bound for $D$ in
Theorem \ref{mythm}, and $\rho_3(D)$, the radius of the ball on which
the convergence in Statement \ref{newit3} is guaranteed, are impossible
to compute without knowing the functions $\rho\mapsto\d(B_\rho),
\rho\mapsto\Delta(B_\rho)$ precisely. Thus Statement \ref{newit4}
gives more easily used, if less sharp, lower bounds on these numbers.
The analogous statement for $\rho_1$ will be used in \S 6 when we
estimate the convergence rate of the sequence $\{\Psi_Q^n(p_0)\}$.
}

\end{remark}

\begin{remark}{\rm 
As $D\to 0$, the  numbers $\rho_3(D)$ and $\rho_4(D)$ increase. Thus, the
smaller the diameter of the set $Q$, the larger the set on which the
theorem shows that the iterates $\Psi_Q^n$ converge, and the larger
the set on which the critical point of $f_Q$ is guaranteed to be
unique.  Also note that $\lim_{D\to 0}\rho_3(D)=\lim_{D\to
0}\rho_4(D)=\sup\{\rho\in [0,r_1]
\mid \k(\rho,0)<1\}$---a considerably larger number than 
$D_{\rm max}=\sup\{\rho\in [0,r_1] \mid \k(\rho,\rho)<1\}$, which is
the upper bound we would have found for the radii of the balls
$B_{\rho_3}, B_{\rho_4}$ in statement \ref{newit3} and
had we not separated the roles of the variables $\rho$
and $D$ in defining $\rho_3$ and $\rho_4$ (i.e. if we had used
``$2\rho$'' in place of ``$\rho+D$'' in (\ref{bddbyh}) and
(\ref{supbound})).  }
\end{remark}

\ni{\bf Proofs of Corollary \ref{weakkarch} and Theorem \ref{mythm}:}
Statements \ref{newit1} and \ref{newit4} of the theorem just restate
some of the conclusions of Lemma \ref{slemma} for easy reference.
Statement \ref{it2} follows from (\ref{kq}), since $s(\rho)>0\iff
\k(\rho)<1$.  Statement \ref{newit3} of
Theorem \ref{mythm} and the existence portion of Corollary
\ref{weakkarch} follow from Theorem \ref{quant_cor} applied to
$U=B_{\rho_4},B=B_\rho,$ since for $\rho_1<\rho<\rho_3$ the fact that
$D < s(\rho)$ ensures that the condition (\ref{kplt1}) is met.  The
conclusion that $\Bar{q}\in\overline{B_{\rho_0}}\cap\ohull(Q)$ just
combines Corollary \ref{primcom} with Karcher's bound
(\ref{karchdist}).

Integrating (\ref{hessbounds}) over $Q$ implies that $\hess(f_Q) > 0$
on $B_\rho$ provided that
$(\rho+D)\cdot\max(0,\tilde{\Delta}(\rho))^{1/2}<\pi/2$, a condition
that Lemma \ref{slemma} (statement \ref{rho4bound}) ensures is met
with $\rho=\tilde{\rho}_4$. Hence Lemma \ref{unimin} implies that any
critical point of $f_Q$ in $B_{\tilde{\rho}_4}$ is unique and
minimizes $f_Q$ on this ball (in fact, on $\regstar(Q)$ by Corollary
\ref{primcom}), proving the remainder of Corollary
\ref{weakkarch}.
\qedns

Theorem \ref{mythm} gives us an algorithm for computing the center of
mass to any desired accuracy: start with some point $q$, and compute
the iterates $\Psi_Q^n(q)$.  As mentioned earlier, when $Q$ is a
finite set of points, it is natural to initialize the algorithm at
some point of $Q$. This motivates the following corollary. In many cases of
interest the ambient manifold is highly symmetric and the quantities
$\rregcvx(q), \tilde{D}_{\rm crit}(q)$ below are independent of $q$,
enabling a much simpler statement of the corollary.

\begin{cor} \label{easycor}
Let $Q\subset M$, $\mu$ a probability measure on $Q$.  For simplicity
let constants $\tilde{\Delta}\ident\Delta(M), \tilde{\d}\ident\d(M)$
be global upper and lower bounds on sectional curvature.  For $q\in Q$
let $D_{q}(Q)=\sup\{d(q,q_1)\mid q_1\in Q\}$, let $\tilde{D}_{\rm
crit}(q)=\tilde{D}_{\rm crit}(q, \rregcvx(q))$ be as in
(\ref{defdcritt}).  If for at least one point $q_0\in Q$ we have
$D_{q_0}(Q)<\tilde{D}_{\rm crit}(q_0)$, then the center of mass
$\Bar{q}$ of $(Q,\mu)$ exists, and equals $\lim_{n\to\infty}
\Psi_Q^n(q)$ for every $q\in Q$.  In particular this conclusion holds
for any $q_0\in Q$ if $\diam(Q)<\tilde{D}_{\rm
crit}(Q):=\inf\{\tilde{D}_{\rm crit}(q_0)\mid q_0\in Q\}$.
\end{cor}

\pf The hypotheses imply that $Q\subset B_D(q_0)$, where $D=D_{q_0}(Q)$.
Letting $\e=\tilde{D}_{\rm crit}(q_0)-D_{q_0}(Q)$ and defining
$\tilde{\rho}_3=\tilde{\rho}_3(q_0,\rregcvx(q_0)-\e/2;D)$ as in
(\ref{defrho3t}), we have $Q\subset B_{\tilde{\rho}_3}(q_0)$ since
$D<\tilde{\rho}_3$.  Hence statement \ref{newit3} of Theorem \ref{mythm}
implies the result.
\qedns

Corollary \ref{easycor_intro} follows immediately.

Centering the underlying convex regular superdisk at a point of $Q$ as
in Corollary \ref{easycor}, while practical, is wasteful in terms of
the restriction on the diameter of $Q$. Any set $Q$ satisfying the
hypotheses of Theorem \ref{mythm} has a {\em (convex regular)
circumradius} $\circumrad(Q)$: the supremum of the radii of open,
strongly convex, regular geodesic balls containing $Q$.  For $\diam(Q)$
sufficiently small (in particular, if $Q$ admits a convex regular
superdisk centered at one of its points) $\circumrad(Q)<\diam(Q)$,
and the conclusion of Corollary \ref{easycor} remains valid if
$\diam(Q)$ is replaced by $\circumrad(Q)$ and if $\tilde{D}_{\rm crit}(Q)$
is replaced by $\tilde{D}_{\rm crit}(p_0)$, where $p_0$ is the
``circumcenter''.  
As a practical matter, the circumcenter is no easier to find
than the center of mass, so that this strengthening of Corollary
\ref{easycor} is only useful if one has a uniform bound on
$\rregcvx(p)$ (and therefore on $\tilde{D}_{\rm crit}(p)$) for $p$ in
an appropriate neighborhood of $Q$.  We will discuss this more
quantitatively in \S \ref{examples}.  


\setcounter{equation}{0}
\section{Rapid convergence of the algorithms}\label{conv_rate}

Given an iterable map $F$, let $\It(F)$ denote the algorithm ``iterate
$F$''.  Under any contracting-mapping algorithm, the sequence of
successive distances from one point to the next converges
geometrically. However, it is well known that Newton's method does
even better; each successive distance is bounded by a constant times
the square of the preceding one.  In this section we examine the
convergence rates of algorithms of the form $\It(\Psi_Y)$ and
$\It(\Phi_X)$ in general (where $\Psi_Y$ and $\Phi_X$ are as in
Theorem \ref{quant_cor}), and of the averaging algorithm
$\It(\Psi_{Y_Q})$ of Theorem \ref{mythm} and Corollary \ref{easycor}
in particular.  We will see that while the convergence rate of
$\It(\Psi_Y)$ for general $Y$ is only geometric (although with a
smaller ratio than $\k(\Psi_Y)$), the algorithms $\It(\Phi_X)$---more
closely related to the flat-space Newton's method---have the same
quadratic behavior as their flat-space cousins.  The averaging
algorithm falls somewhere in between: we obtain only geometric
convergence,  but with a very small ratio, provided that
$\diam(Q)$ is small enough.

Throughout this section, notation will be as in Theorem
\ref{quant_cor}. We denote the sequence of iterates $\{\Psi_Y^n(p_0)\}$ or
$\{\Phi_X^n(p_0)\}$ by $\{p_n\}$.  For any algorithm of the form
$\It(\Psi_Y)$, the following proposition shows that the rate at which
$d(p_n,p_{n+1})\to 0$ is completely controlled by bounds on $\na
Y+I$. 

\begin{prop} \label{control_lemma}
Let $U$ be a convex set preserved by $\Psi_Y$, let $p_0\in U$, and for
$n>0$ let $p_n=\Psi_Y^n(p_0)$. Then
\be\label{control}
d(p_{n+1},p_n) \leq (\sup_{p\in U}\|(\na Y +I)_p\|)\ d(p_n,p_{n-1}).
\ee
\end{prop}

\pf
From the definition of $\Psi_Y$, we have
\be\label{yn}
d(p_{n+1},p_n)=\|Y_n\|.
\ee
To analyze how $\|Y_n\|$ changes when we increment $n$, fix $n$ and
let $\g:[0,1]\to M$ be the geodesic from $p_{n}$ to $p_{n+1}$ with
initial velocity $Y_{n}$; thus $p_{n+1}=\g(1)$, $Y_n=Y_{\g(0)}$, and
$Y_{n+1}=Y_{\g(1)}$.  Let $\P_{\g(t)\to\g(0)}$ denote the operator of
parallel transport along $\g$, with direction reversed, from $\g(t)$
back to $\g(0)$, let $A_p=(\na Y+I)|_p\in {\rm End}(T_pM)$, and let
$\e_1=\sup_{p\in U}\|A_p\|$.  Then
\bestar
\frac{d}{dt}(\P_{\g(t)\to\g(0)} (Y_{\g(t)}) +t Y_{\g(0)})
&=&\P_{\g(t)\to\g(0)}(\na_{\g'(t)}Y) +\g'(0)\\
&=&\P_{\g(t)\to\g(0)}(A_{\g(t)}(\g'(t))),
\eestar
since
$\g'$ is parallel along $\g$, and hence
\be
\P_{\g(t)\to\g(0)} (Y_{\g(t)}) +(t-1) Y_{\g(0)}=\int_0^t
\P_{\g(t_1)\to\g(0)}(A_{\g(t_1)}(\g'(t_1))) dt_1.
\ee
The integrand is bounded in norm by $\|A_{\g(t_1)}\| \|\g'(t_1)\|
=\|A_{\g(t_1)}\| \|Y_n\|$. Hence
\bearray\label{bd1}
\|Y_{\g(t)}\|= \|\P_{\g(t)\to\g(0)} (Y_{\g(t)})\| &\leq&
(1-t +\int_0^t\|A_{\g(t_1)}\|\ dt_1)\ \|Y_n\| \\
&\leq& (1-t+\e_1 t)\|Y_n\|.\label{bd2}
\eearray
Inserting $t=1$ we find
$\|Y_{n+1}\|\leq \e_1 \|Y_n\|$, 
and hence 
\be
\label{dbound1}
d(p_{n+1},p_{n})\leq \e_1 d(p_n,p_{n-1}).
\ee
\qedns

Thus in algorithms of the form $\It(\Psi_Y)$, successive distances
decrease geometrically, but with ratio $\e_1$---a number smaller than
the contraction constant $\k(\Psi_Y)$ in (\ref{defkpsiy}), and one
whose only dependence on curvature is through $Y$ itself. 

To analyze the algorithms $\It(\Phi_X)$, proceed as above but with
$Y=-(\na X)^{-1}X$; continue writing $A=\na Y+I$.
In this case, for $p\in U$ and $v\in T_pM$, from
(\ref{nay}) we have $A_p(v)=B_p(v)(Y_p)$, where $B_p(v)=-((\na
X)^{-1}\circ (\na_v\na X))|_p$.  Thus, pointwise we have
\be\label{k3}
\|A\|\leq k_3\|Y\|
\ee
where $k_3=k_1^{-1}k_2$. Inserting this bound into (\ref{bd1}) with
$t=1$, and
using (\ref{bd2}) in the new integrand, we obtain
$\|Y_{n+1}\|\leq \frac{1}{2}k_3(\e_1+1) \|Y_n\|^2$
where now $\e_1=k_1^{-1}\e$.  Thus, with $k_4=k_3(\e_1+1)/2$, we have
\be
d(p_{n+2},p_{n+1})\leq k_4 d(p_{n+1},p_n)^2,
\ee
the same quadratic falloff as in flat-space Newton's method.

Note that the preceding analysis applies to any algorithm for which
(\ref{k3}) holds, a condition intermediate between Case 1 and Case 2
of Theorem \ref{quant_cor}.

The convergence rates of $\It(\Psi_Y)$ and $\It(\Phi_X)$ can also be
compared as follows. With the constants as named above, assume that
for $\Psi_Y$ that$\e_1<1$, and for $\Phi_X$ that $k_4\e_1<1$.
Then for the algorithm $\It(\Psi_Y)$, we have
\be\label{dbound2}
d(p_{n+1},p_n)\leq d(p_1,p_0)\e_1^{n}<\e_1^{n+1},
\ee
whereas for $\It(\Phi_X)$ we have
\be\label{bn}\label{dbound3}
d(p_{n+1},p_n)\leq k_4^{-1}(k_4 d(p_1,p_0))^{2^n}<k_4^{-1}(k_4
\e_1)^{2^n}
\ee
(if $k_4=0$, interpret (\ref{dbound3}) as $d(p_{n+1},p_n)=0$.)

In the proof of the Contracting Mapping Theorem (Theorem \ref{cmt}), to
obtain convergence of the sequence $\{p_n=F^n(p_0)\}$, it suffices to know
that (i) $d(p_n,p_{n+1})\leq \k d(p_{n-1},p_n)$ for all $n\geq 1$, and
(ii) $d(p_0,p_1)<(1-\k)\rho$.  One does not need to know that $F$
is a contraction on the whole ball $B$ unless one wants to prove
uniqueness of the fixed point and convergence of the sequence with
other starting points.  Thus the analysis above leads immediately to
the following existence/convergence theorem to supplement Theorem
\ref{quant_cor}. 

\begin{thm}\label{strongerthm1}
Let $B=B_{\rho}(p_0)\subset M$ be a convex ball. Assume either of the
sets of hypotheses listed as ``Case 1'' and ``Case 2'' in Theorem
\ref{quant_cor}, with $U$ replaced by the ball $B$. 
In Case 1, let $F=\Psi_Y$; in Case 2 let $F=\Phi_X$.
Assume in addition the following:

\ss{\em Case 1.} Assume $\|Y(p_0)\|<(1-\e_1)\rho$.

\ss{\em Case 2.} Let $k_4=k_1^{-1}k_2(k_1^{-1}\e+1)/2$ and if $k_4\neq
0$ assume that
\ben
\sum_{n=0}^\infty k_4^{-1}(k_4k_1^{-1}\|X(p_0)\|)^{2^n} <\rho.
\een

\ssn
Then in each case the sequence $\{F^n(p_0)\}$ lies in $B$ and
converges to a fixed point of $F$ that lies in $B$.

The distances $d(p_{n+1},p_n)$ in Case 1 have the exponential falloff
given by (\ref{dbound2}), and in Case 2 have the super-exponential
falloff given by (\ref{dbound3}). \qedns
\end{thm}

Theorem \ref{strongerthm1} is most useful when one knows ahead of time
that there is at most one fixed point.  This is exactly the case for
averaging algorithm $\It(\Psi_{Y_Q})$ used in \S \ref{construct},
since we do not need the contracting-mapping apparatus to prove
uniqueness---given existence, we already know from Kendall's theorem
that if $Q$ is contained in regular geodesic ball $B$ then
$\Psi_Q:=\Psi_{Y_Q}$ has at most one fixed point in $B$.  This leads
immediately to the following strengthening of certain portions of
Theorem \ref{mythm}.

\begin{thm} \label{strongerthm2} 
Let $p_0\in M$, $0<r_1\leq\rregcvx(p_0)$; for $0<\rho\leq r_1$ write
$B_\rho$ for $B_\rho(p_0)$. Let $\tilde{\Delta}(\cdot),
\tilde{\d}(\cdot)$ be continuous monotone upper and lower bounds on
curvature as in Definition \ref{def44}(b). Define numbers
$\tilde{D}'_{\rm crit}, \tilde{D}'_{\rm max}, \tilde{\rho}'_{\rm
crit}$ and $\tilde{\rho}'_{i}$ analogously to the numbers defined in
Lemma \ref{slemma}, but with $\tilde{s}$ replaced by the function
\be\label{defsseq}
\tilde{s}_{\rm seq}(\tilde{\Delta},\tilde{\d};\rho,D)=
(1-\tilde{\k}_{\rm seq}(\tilde{\Delta},\tilde{\d};\rho,D))\rho
\ee
where 
\be\label{defkseq}
\tilde{\k}_{\rm seq}(\tilde{\Delta},\tilde{\d};\rho,D))
=\psimax(\tilde{\d}(\rho),\tilde{\Delta}(\rho),\rho+D).
\ee
Then Statements \ref{newit1} and \ref{newit4} 
of Theorem \ref{mythm} hold with
$\tilde{D}_{\rm crit}, {D}_{\rm crit}, \tilde{\rho}_{i}$, and
${\rho}_{i}$ replaced by \newline $\tilde{D}'_{\rm crit}, {D}'_{\rm
crit},
\tilde{\rho}'_{i}$, and ${\rho}'_{i}$
respectively.  Assume that $Q\subset \Bar{B_D}$ (or more generally
$\|Y_Q(p_0)\|\leq D$) and that either

\ss (i) $D<\tilde{D}'_{\rm crit}$, or 

\ss (ii) $D<\tilde{D}'_{{\rm max}}$ and  $\Psi_Q$ is tethered to $Q$
(see Definition \ref{deftether}).

\ss 
Then the sequence of iterates $\{\Psi_Q^n(p_0)\}$ converges to the primary
center of mass $\Bar{q}$ of $(Q,\mu)$, and $\Bar{q}$ lies in
$\Bar{q}\in\overline{B_{\rho_0(D)}}\cap\ohull(Q)$. The entire sequence
lies in $\Bar{B_{\rho'_{1}(D)}}$ (hence in the $D$-independent ball
$B_{\rho'_{\rm crit}}$) if (i) holds, and in $B_{{\rho}'_{4}(D)}$ if
(ii) holds. \qedns
\end{thm}

We have a corresponding strengthening of Corollary \ref{easycor}:

\begin{cor}\label{strongerthm3}
Corollary \ref{easycor} remains true if the numbers $\tilde{D}_{\rm
crit}(q)$ are replaced by the larger numbers $\tilde{D}'_{\rm
crit}(q)$ defined in Theorem \ref{strongerthm2}. \qedns
\end{cor}

\ss
For the map $\Psi_Q=\Psi_{Y_Q}$ 
used in Theorem \ref{strongerthm2} and Corollary
\ref{strongerthm3} we have a bound on
the endomorphism $A$ that, while not as strong as (\ref{k3}), is
better than for the general $\Psi_Y$.  From (\ref{bddbyh}) and
(\ref{H2}), 
if $Q\subset B_D(p_0)$ then on
$B_\rho(p_0)$ we have
\be\label{abound}
\|A\| \leq \psimax(\d(B_\rho(p_0)),\Delta(B_\rho(p_0)),\rho+D) 
=\frac{1}{3}|K| (\rho+D)^2+ O(|K|^2 (\rho+D)^4)
\ee
where $|K|=\max(\d(B_\rho(p_0)),\Delta(B_\rho(p_0))$.
Initialize the algorithm at a point $p_0\in Q$ as in Corollary
\ref{easycor}, let $D=\diam(Q)$, and assume that
$D<D_{\rm crit}(Q)$ as in the corollary.  From Theorem \ref{strongerthm1}
$\Psi_Q$ preserves the convex ball $\Bar{B_{\rho_1}(p_0)}$, where $\rho_1(D)$
is the smallest positive number $\rho$ satisfying $s(\rho,D)=D$, and
hence when applying the bound (\ref{abound}) in the analysis of
$\{\Psi_Q^n(p_0)\}$ it suffices to take $\rho=\rho_1(D)$.  Since
$s(\rho,D)=\rho(1-O(|K|(\rho+D)^2))$, for $D$ small we have
$\rho_1(D)=D(1+O(|K|D^2))$.  Thus
\be\label{abound2}
\|A\| \leq 
 \frac{4}{3}|K|\diam(Q)^2 +O(|K|^2\diam(Q)^4),
\ee
which we can use for $\e_1$ in (\ref{dbound1}) and (\ref{dbound2}). 
Thus for any $\e_2>0$, if $|K|\cdot\diam(Q)^2$ is small enough we have
\bearray
\label{dbound4}
d(p_{n+1},p_{n})&\leq& (\frac{4}{3}+\e_2)|K|\diam(Q)^2
d(p_n,p_{n-1})\\
&=&k_5\ \diam(Q)^2 d(p_n,p_{n-1}), 
\eearray
so in place of (\ref{dbound2}) we can write
\be\label{dbound5}
\frac{d(p_{n+1},p_n)}{\diam(Q)^n}\leq d(p_1,p_0)(k_5\diam(Q))^n.
\ee

In other words, as $\diam(Q)\to 0$, the falloff rate of successive
distances in the averaging algorithm is geometric even relative to
$\diam(Q)$.  The bound (\ref{dbound3}) shows that we would get even
faster convergence to the center of mass if we iterated the map
$\Phi_{Y_Q}$ instead of $\Psi_{Y_Q}$. However, as a practical tool
$\Phi_{Y_Q}$ has the disadvantage that one must compute and invert
$\na Y_Q$, which may be difficult even if $M$ has constant curvature,
whereas for many more general spaces the algorithm $\It(\Psi_{Y_Q})$
is easily programmable.

\begin{remark}{\rm
Since $A=\na Y +I$, for $\diam(Q)$ small we can think of
(\ref{abound}) as asserting that the vector field $Y_Q$ is, in some
sense, very nearly linear. From this point of view it is no surprise
that the convergence of the algorithm is so rapid---what we are using
is almost Newton's method for an almost linear function.
}
\end{remark}

As $D\to 0$, the bound (\ref{dbound4}) can be improved by using the
circumradius of $Q$ instead of its diameter in this estimate (see the
discussion after Corollary \ref{easycor}).  In
$\bfr^n$, one always has $\circumrad(Q)\leq
\sqrt{\frac{n}{2(n+1)}}\ \diam(Q)$, with a regular $n$-simplex an
extremal configuration.  In a general Riemannian manifold, if we
restrict attention to sets $Q$ contained in a subset $U$ on which the
there are bounds on the curvature and a positive lower bound on the
injectivity radius, then as $D\to 0$ the
number $\sup\{\circumrad(Q)/\diam(Q)\mid Q\subset U,\ 0<\diam(Q)\leq
D\}$ tends to its Euclidean value. Thus we obtain an asymptotic bound
$\e_1\sim \frac{2}{3}\frac{n}{n+1} \Delta D^2$, where $n=\dim(M)$.


\setcounter{equation}{0}
\section{Averaging in the case of non-negative curvature}
\label{examples}

When $(M,g)$ has curvature of a fixed sign, the definitions of the
critical radii in Theorem \ref{strongerthm2} and Corollary
\ref{strongerthm3} simplify, since we can globally replace 
$\psimax(\d_{p,\rho},\Delta_{p,\rho},\rho+D)$ in (\ref{defkseq}) by
either either $\psi(\tilde{\Delta}(\rho),\rho+D )$ or
$\psi(\tilde{\d}(\rho),\rho+D)$.  In this section we assume that the
curvature is non-negative, which is true in all shape spaces and
size-and-shape spaces.

The goal of this section is to estimate the critical radii appearing
in Theorem \ref{strongerthm2} as well as the convergence rate of the
averaging algorithm (not merely the asymptotics of this rate as
$\diam(Q)\to 0$). To simplify the estimates further, we will assume a
uniform upper bound $\tilde{\Delta}\ident\Delta$ on sectional
curvature in all the balls that appear in this section, and a uniform
lower bound $r_1$ on the regular convexity radius of the center of any
such ball.  We assume $\Delta>0$ strictly since the flat case is not
very interesting, the algorithm converging at the first iteration.

Notation in this section will be for the most part as in \S\S
\ref{construct}--\ref{conv_rate}, but it is convenient
to define rescaled variables $\bar{\rho}=\Delta^{1/2}\rho,
\bar{D}=\Delta^{1/2}D$, and a rescaled function
$\bar{s}=\Delta^{1/2}\tilde{s}$ of the rescaled variables (where in
the definition of $\tilde{s}$ we take $\tilde{\d}\ident 0,
\tilde{\Delta}\ident\Delta$). We also write $\bar{\k}$ for
$\tilde{\k}$ expressed in terms of the rescaled variables.
We suppress all the parameters except $D$ and $\rho$ in most
formulas below. 

Fix $p_0\in M$ and let $x=\bar{\rho}+\bar{D}$. Then
\be\label{nonneg1}
\bar{\k}(\bar{\rho},\bar{D})
=\hat{\k}(x):=
\psi(1,x)= 1-x\cot x = \frac{1}{3}x^2+O(x^4)
\ee
and
\be
\bar{s}(\bar{\rho},\bar{D})=(1-\bar{\k}(\bar{\rho},\bar{D}))\bar{\rho}.
\ee
Since $\tilde{\Delta},\tilde{\d}$ are constant, $\bar{s}$ is
differentiable, so the rescaled pair $(\bar{\rho}_{\rm crit},
\bar{D}_{\rm crit})$ from Lemma
\ref{slemma} can be
characterized as the unique solution of the system of equations
\bearray
\bar{s}(\bar{\rho},\bar{D})&=&\bar{D},\\
\frac{\partial\bar{s}}{\partial\bar{\rho}}(\bar{\rho},\bar{D}) &=& 0
\eearray
in $(0,\pi/2)\times (0,\pi/2),$ provided that $\bar{\rho}_{\rm
crit}$ as defined this way is less than $\Delta^{1/2}r_1$.
For this system of equations, Maple's {\tt fsolve} routine\footnote{
All numerical calculations in this section were done with Maple.}
yields $
\bar{\rho}'_{\rm crit}\approx .6816\gtrsim .2169\pi,
\bar{D}'_{\rm crit}\approx .3952 \gtrsim .1258\pi.$
Thus
\be\label{calc1}
\tilde{\rho}'_{\rm crit}\geq \min(r_1,.2169\pi \Delta^{-1/2}), \ \ 
\tilde{D}'_{\rm crit}\geq \min(r_1, .1258\pi \Delta^{-1/2}).
\ee
From these numbers we also compute
$\tilde{\rho}'_4(\tilde{D}'_{\rm crit})\approx \min(r_1,1.1566 \Delta^{-1/2})
\approx .3682 \pi\Delta^{-1/2}.$
Centering all balls below at $p_0$ and writing $B_\rho$ for
$B_\rho(p_0)$, we recall what the numbers just computed tell us: from
Theorem \ref{strongerthm2}, for any $(Q,\mu)$ with $\Bar{Q}$ in the
ball of radius $\tilde{D}'_{\rm crit}$, and any $p$ in the ball of
radius $\tilde{\rho}'_{\rm crit}$, the sequence $\{\Psi_Q^n(p)\}$
converges to the primary center of mass of $(Q,\mu)$.  If $\Psi_Q$ is
tethered to $Q$, to conclude convergence we need only assume that
$\Bar{Q}$ and $p$ lies in the balls of radius $\tilde{D}'_{\rm crit}$
and $\tilde{\rho}'_{\rm 4}(\tilde{D}'_{\rm crit})$ respectively.

If $Q\subset B_D$ then as $D\to 0$, the algorithm converges on larger
and larger sets, the balls of radius $\tilde{\rho}'_3(D)$ (or
$\tilde{\rho}'_4(D)$ in the tethered case). These radii approach
$\tilde{\rho}'_3(0)=\tilde{\rho}'_4(0)=\min(r_1,(\pi/2)\Delta^{1/2})$.
Thus as $D\to 0$ we get convergence on balls of radius arbitrarily
close to (but smaller than) the largest radius for which Kendall's
theorem (Theorem \ref{kendthm}) guarantees uniqueness of the center of
mass.

\begin{remark}{\rm
Corollary \ref{weakkarch}, the existence/uniqueness theorem given by
the contracting-mapping approach, guarantees existence of the center
of mass of a distribution supported in a ball of radius
$\tilde{D}'_{\rm crit}$; in Karcher's result, the $.1258\pi$ in
(\ref{calc1}) is replaced by the better $\pi/4$. To compare the
uniqueness statement in Corollary \ref{weakkarch} with those of
Karcher and Kendall, we cannot use the radii above, coming from
Theorem \ref{strongerthm2}, but must go back to those in Theorem
\ref{mythm}. This has the effect of replacing $\psi(1,x)$ in
(\ref{nonneg1}) by $\phi_-(x)+\psi(1,x) = \frac{2}{3}x^2 +O(x^4)$. 
In this case we analogously compute $\tilde{D}_{\rm crit}\approx
\min(r_1, .0904\pi \Delta^{-1/2})$ and $\tilde{\rho}_4(\tilde{D}_{\rm
crit})\approx \min(r_1, .2777\pi \Delta^{-1/2})$\footnote{These
numbers increase slightly if $(M,g)$ is further assumed to be locally
symmetric, since instead of $\phi_-(x)$ we can then use the smaller
quantity $\phi_+(x)=\phi_-(x)-\frac{1}{15}x^4+O(x^6)$.  In this case
we can replace $.0904\pi$ by $.0932\pi$, and $.2777\pi$ by $.2991\pi$.
The improvement is so marginal because $\phi_+(x)$ and $\phi_-(x)$
differ by only $\frac{1}{15}x^4+O(x^6)$.}.  Corollary \ref{weakkarch}
implies that if $\Bar{Q}$ is contained in the ball of radius
$\tilde{D}_{\rm crit}$, then $(Q,\mu)$ has a unique center of mass in
the ball of radius $\tilde{\rho}_4(\tilde{D}_{\rm crit})$.  Thus in
the non-negative curvature case, for $D<\tilde{D}_{\rm crit}$ the
contracting-mapping approach, while giving not as strong a uniqueness
statement as in Kendall's theorem, gives a slightly stronger
statement than in Karcher's original theorem, which has
only $\pi/4$ in place of our worst-case constant $.2777\pi$.
}
\end{remark}

We next estimate the convergence rate of $\{p_n=\Psi_Q^n(p_0)\}$,
assuming that $\Bar{Q}$ lies in the ball of radius ${\tilde{D}'_{\rm
crit}}$.  From Theorem
\ref{strongerthm2} the sequence stays in the ball 
of radius ${\rho'_{\rm crit}(D)}$, on which, letting $A=\na Y_Q +I$
and writing $x_{\rm crit}={\bar{\rho}'_{\rm crit}(D)}+{\bar{D}'_{\rm
crit}}$, the bound (\ref{abound}) gives
\be
\|A\|\leq 
\psimax(0,1,x_{\rm crit})=\hat{\k}(x_{\rm crit})
\lesssim .4202
\ee
Hence we obtain the geometric convergence rate (\ref{dbound1}) with
$\e_1=.4202$. 

If we start with $p_0\in Q$ and assume $D=\diam(Q)<\tilde{D}'_{\rm crit}$ as
in Corollary \ref{strongerthm3}, then as $D$ decreases we can sharpen the
convergence-rate estimate by replacing $x_{\rm crit}$ with
$\bar{\rho}'_1(D)+\bar{D}$ in the previous estimate.  Since
$\hat{\k}$ is monotone increasing on $[0,x_{\rm crit}]$, and
$\tilde{s}(\rho_1(D),D)=D$, we have $\rho_1(D)\leq
\frac{D}{1-\hat{\k}(x_{\rm crit})}:=c_1D \leq 1.725 D$. 
The function $x\mapsto \hat{\k}(x)/x^2$ is monotone increasing on
$[0,\pi)$, so for $x\in [0,x_{\rm crit}]$ we have $0\leq \hat{\k}(x)\leq
(\hat{\k}(x_{\rm crit})/x_{\rm crit}^2)x^2 :=c_2x^2$.  Thus
$\|A\| \leq c_2(1+c_1)^2\Delta D^2\leq 2.690 \Delta D^2,$ so we can
take $\e_1=2.690 \Delta D^2$ in (\ref{dbound1}) and (\ref{dbound2}).

As $D\to 0$, this can be improved further---(\ref{dbound4}) gives a
bound on $\e_1$ asymptotic to $\frac{4}{3}\Delta D^2$, and as noted at
the end of \S \ref{conv_rate} this can even be reduced to 
$\frac{2}{3}\frac{n}{n+1} \Delta
D^2$, where $n=\dim(M)$.

Finally, we consider two simple examples: round spheres and complex
projective spaces, with standard metrics.  If $M$ is a round sphere of
radius $R$, then the curvature is constant and equal to $R^{-2}$, and
$\rcvx(M)=\rreg(M)=\pi R/2$. Hence  we can take
$\Delta^{-1/2}=R$ and erase ``min'',``$r_1$'' and the tildes in all
the estimates above; e.g. in place of (\ref{calc1}) we have simply
\be\label{sphere}
{\rho}'_{\rm crit}\geq .2169\pi R, \ \ 
{D}'_{\rm crit}\geq .1258\pi R
\ee
Similarly, ${\bf C}P^n$ with a Fubini-Study metric (unique up to
scale) is a symmetric space of positive curvature. If we fix the scale
by taking the metric to be the one for which the standard projection
from the unit sphere $S^{2n+1}\to {\bf C}P^n$ is a Riemannian
submersion, then the sectional curvatures of ${\bf C}P^n$ run between
$\d=1$ and $\Delta=4$ if $n\geq 2$ (the curvature is identically 4 if
$n=1$; ${\bf C}P^1$ with this metric is a round sphere of radius
$1/2$). In this case we have $r_{\rm cvx}(M)=\pi/4$ and
$\Delta^{-1/2}=1/2$, so the critical radii are exactly half those for
the unit sphere; bounds are given by (\ref{sphere}) with $R=1/2$.  It
is not hard to show that $\Sigma^k_2$, the shape space of $k$ points
in $\bfr^2$, is exactly $\bfc P^{k-2}$ with this metric (if $k>2$)
\cite{dkendall}, so the numbers above directly relate to the behavior
of the Riemannian averaging algorithm on this shape space.


\setcounter{equation}{0}
\section{Appendix}\label{appendix}

\subsection{Proof and discussion of Proposition \ref{shortgeods}}
\label{app2}

In this subsection, hypotheses and notation are as in Proposition
\ref{shortgeods}. We first prove (\ref{jbound}) and then discuss how
to sharpen this bound for locally symmetric spaces; the bound
(\ref{jbound_symplus}) follows as a special case of this discussion.

\bs
{\bf Proof of (\ref{jbound}).}  $J^\parallel$ and $J^\perp$, the
components of $\hat{J}_v$ parallel and perpendicular to $\g'$, are
themselves Jacobi fields, with $J^\parallel(t)=(at+c)\g'(t)$ for some
$a,c\in\bfr$. Each of $J^\parallel$ and $J^\perp$ satisfies
antidiagonal initial conditions.  In particular, $c=-a$, so
$J^\parallel(1)=0$. Hence $\hat{J}_v(1) = J^\perp(1)$, so it suffices
to prove (\ref{jbound}) under the assumption that $v\perp\g'(0)$,
which we make henceforth.

Let $\{e_i\}_0^{n-1}$, where $n=\dim(M)$, be an orthonormal basis of
$T_pM$ with $e_0=\g'(0)/\|\g'(0)\|$, and extend each $e_i$ along $\g$
by parallel translation. Write $J(t)=\sum_{i=1}^{n-1} f^i(t)e_i(t)$
and let $f:[0,1]\to\bfr^{n-1}$ be the vector-valued function whose
components are the $f^i$; note that $\|f(t)\|_{\rm
Euclidean}=\|J(t)\|$. Then (\ref{jaceq}) simply becomes
\be\label{jaceq2}
f''(t)=A(t)f(t)
\ee
for a certain $(n-1)\times (n-1)$ matrix-valued function $A$ whose
operator norm satisfies
$\|A(t)\|\leq \aku(\g(t))\|\g'(t)\|^2$.  The norm of $\g'(t)$ is constant
and equal to the length $r$ of $\g$.  Letting $b=\aku(\g)$, we
therefore have $\|A(t)\|\leq br^2.$

For $v\in T_pM$ write $v=\sum v^ie_i$, and let $\bar{v}\in\bfr^n$ be
the vector whose components in the standard basis are the $v^i$. The
initial conditions for $\hat{J}_v$ then become
$f(0)=-f'(0)=\bar{v}.$ The unique solution of (\ref{jaceq2}) with
these initial conditions is given explicitly by the series

\bearray\nonumber
f(t) &=& (1-t)\bar{v} +\int_0^t\int_0^{t_2} (1-t_1)A(t_1)\bar{v}\ dt_1dt_2 
\\ \nonumber && +\int_0^t\int_0^{t_4}\int_0^{t_3}\int_0^{t_2}
(1-t_1)A(t_3)A(t_1)\bar{v}\ dt_1dt_2dt_3dt_4 \\ 
\nonumber
&& +\dots + \int\dots\int_{0\leq t_1\leq t_2 \dots \leq
t_{2m}\leq t}
(1-t_1)A(t_{2m-1})A(t_{2m-3})\dots A(t_1)\bar{v}\ dt_1\dots
dt_{2m}\\
\label{series} &&  +\dots
\eearray
(This series converges in norm uniformly on any compact $t$-interval.)
In the $2m$-fold integral the integrand is bounded in norm by
$(1-t_1)b^mr^{2m}\|v\|$ provided $0\leq t\leq 1$, the only case we are
interested in.  Integrating explicitly, we obtain
$\|v\|b^mr^{2m}(\frac{t^{2m}}{(2m)!} -
\frac{t^{2m+1}}{(2m+1)!})$ as an upper bound on the $2m$-fold
integral. Hence for $0\leq t\leq 1$ we have
\bestar
\|f(t)\|
&\leq& \sum_{m=0}^\infty b^mr^{2m}(\frac{t^{2m}}{(2m)!} - 
\frac{t^{2m+1}}{(2m+1)!})\|v\| \\
&=&(\cosh (b^{1/2}rt) -\frac{\sinh(b^{1/2}rt)}{b^{1/2}r})\|v\|.
\eestar
Plugging in $t=1$, the bound (\ref{jbound}) follows.
\qed

In contrast to more frequently-seen bounds on Jacobi fields, the sign
of the sectional curvature does not play a role in (\ref{jbound}). The
reason is the anti-diagonal initial condition, which in Euclidean
space leads to $J(1)=0$. If $M$ is positively curved, then $\|J\|$ can
reach 0 before time 1 and then grow again, so that $\|J(1)\|$ cannot
be bounded by its Euclidean analog.  However, while it is not obvious
how to get the best bound in Proposition \ref{shortgeods} for general
manifolds, or even for nonnegatively curved manifolds, the analysis
simplifies considerably for locally symmmetric spaces (manifolds whose
Riemann tensor is covariantly constant; examples are $S^n$ and ${\bf
C}P^n$). In this case the matrix $A(t)$ in
(\ref{jaceq2}) is a constant symmetric matrix $r^2\hat{A}$, and the
solution (\ref{series}) collapses to
\be\label{collapse}
f(t)=(\ch(t^2r^2\hat{A})-t\ \sh(t^2r^2\hat{A}))\bar{v}
\ee
(see Table 1 in \S 2.)  Hence in this case (\ref{jbound}) can be
improved to
\be \label{jbound3}
\| \hat{J}_v(1) \| \leq \|\ch(r^2\hat{A})-\sh(r^2\hat{A})\|\ \|v^\perp\|.
\ee
We can always choose an orthonormal basis in which the matrix
$\hat{A}$ in the proof above is diagonal, say $\hat{A}={\rm
diag}(\l_1,\dots,\l_{n-1})$. Then $\ch(r^2\hat{A})-\sh(r^2\hat{A})$
becomes a diagonal matrix with entries ${\rm
sign}(\l_i)\cdot\phi_{{\rm sign}(\l_i)}(|\l_i|^{1/2}r).$ The sectional
curvatures of $M$ range between $\d\leq\min\{{\l_i}\}$ and $\Delta\geq
\max\{\l_i\}$ (we would have equality here if we replaced $\d$ and
$\Delta$ by the minimum and maximum sectional curvatures achieved on
2-planes tangent to $\g$) and $\phi_\pm$ are increasing functions on
appropriate intervals: $\phi_-$ on $[0,\infty)$ (the Taylor
coefficients are all nonnegative), $\phi_+$ on $[0,x_0]$, where
$x_0\approx 0.87 \pi$ is the first positive solution of $(x^2-1)\sin
x+ x\cos x=0$.  Hence
\be\label{replphi}
\|\ch(r^2\hat{A})-\sh(r^2\hat{A})\|\leq
\left\{ \begin{array}{ll} 
\phi_+(\Delta^{1/2}r) & {\rm if}\
0\leq\d\leq \Delta\
{\rm and}\ \Delta^{1/2}r\leq x_0, \\
\max(\phi_-(|\d|^{1/2}r), \phi_+(\Delta^{1/2}r)) & {\rm if}\
\delta\leq 0<\Delta\
{\rm and}\ \Delta^{1/2}r\leq x_0, \\
\phi_-(|\d|^{1/2}r) &  {\rm if}\
\delta\leq\Delta<0. \end{array}\right.
\ee

Thus for a locally symmetric space we can replace
$\phi_-(r\aku(\g)^{1/2})$ in (\ref{jbound}) by the appropriate line of
(\ref{replphi}); the top line yields (\ref{jbound_symplus}), since
$x_0>3\pi/4$.  
(We chose $3\pi/4$ in Proposition \ref{shortgeods} for simplicity. 
Values of $\phi_+$
that equal or exceed 1 are irrelevant to us since in Theorem
\ref{quant_cor}(b) they lead to a useless bound on $\k$.  The first
positive $x$ for which $\phi_+(x)=1$ is approximately $.74\pi$, so the
restriction $\Delta^{1/2}r\leq 3\pi/4$ more than suffices for our
considerations.)

If $M$ has {\em constant} curvature---i.e. all sectional curvatures
are equal, say to $\Delta$---then the matrix in (\ref{collapse}) is a
multiple of the identity, leading us to sharp equality.  In this
case $\hat{A}=-\Delta I$ so we obtain
\be \label{jbound4}
\| \hat{J}_v(1) \| = \phi_\pm(|\Delta|^{1/2}r)\|v^\perp\|
\ee
where $\phi_+$ is used if $\Delta\geq 0$, and $\phi_-$ if $\Delta<0$.

\subsection{The Hessian of the squared distance
function}\label{hess_app}

Good references for the material in this subsection are
\cite{hild}, \S 5 and \cite{karcher}, Appendix C.

The lemma below was used in Lemma \ref{hesslemma} and Corollary
\ref{weakkarch}. The useful bound (\ref{hessbounds_app}) is
essentially proven in \cite{hild} Chapters 4-5, but is not explicitly
stated in this form. (Theorem 5.2 of \cite{hild} asserts an inequality
that looks identical to (\ref{hessbounds_app}), but because
Hildebrandt's goal in \cite{hild} is a simple upper bound that applies
to {\em all} vectors, not just those orthogonal to $\g'$, he imposes
the requirement $\d\leq 0$.) The block-diagonal decomposition of the
Hessian indicated in the lemma must generally be used in order to get
the sharpest estimates on $\|\na Y+I\|$ when $Y$ is the gradient of a
function of the form $p\mapsto \int_Q f(d(p,q))\ d\mu(q)$.

\begin{lemma}\label{hesslemma_app}
Let $p,q\in M$ with $d(p,q)<\rinj(q)$ and let
$H=\hess(\frac{1}{2}r_q^2)|_p$.  Let $\g:[0,1]\to M$ be the minimal
geodesic from $q$ to $p$, let $u$ a unit vector tangent to $\g$ at
$p$, and let $V_p^\perp\subset T_pM$ be the orthogonal complement of
$\span(u)$.  Let $\d$ and $\Delta$ be lower and upper
bounds, respectively, for the sectional curvatures of $M$ along $\g$;
if  $\Delta>0$ also assume $d(p,q)<\pi\Delta^{-1/2}$.
Then for all $v\in V_p^\perp$ we have the following:
\bearray
H(u,u)&=&1, \label{hesslemma1} \\
H(u,v)&=&0, \label{hesslemma2} \\
h(\Delta,d(p,q))\|v\|^2 &\leq& H(v,v)\leq
h(\d,d(p,q))\|v\|^2.
\label{hessbounds_app} 
\eearray
\end{lemma}

\pf Recall that for any function $f$, vectors $X,Y\in T_pM$, and
an arbitrary smooth extensions of $X,Y$ to vector fields on a
neighborhood of $p$, the covariant Hessian $H_f$ is given by
\be\label{hessdef}
H_f(X,Z)=X(Z(f))-(\na_XZ)(f).
\ee
Let $f=\frac{1}{2}r_q^2$, let $X$ be an extension of the unit tangent
vector field $\g'/\|\g'\|$ and let $Z$ be an extension of $v\in
V^\perp_p$ that is parallel along $\g$. Then (\ref{hesslemma1}) is
trivial, and, since the Gauss lemma implies $Z(r_q)\ident 0$ along
$\g$, (\ref{hesslemma2}) is trivial as well. The bound
(\ref{hessbounds}) can derived from the normal-Jacobi-field estimate
\cite{hild} Theorem 4.2, followed by rescaling the arclength parameter as at
the bottom of \cite{hild} p. 53, and then restricting the proof of
\cite{hild} Theorem 5.2 to the case of vectors orthogonal to the
geodesic.
\qed


\end{document}